\documentclass[a4paper,12pt,mathcal]{article}

\usepackage{amsfonts,eucal,amssymb,amsmath,amsthm,amscd} \input{epsf}

\newtheorem{theoreme}{Theorem}

\newtheorem{corollaire}{Corollary}

\newtheorem{lemme}{Lemma}
\newtheorem{definition}{Definition}

\newtheorem{remarque}{Remark}

\newtheorem{exemple}{Example}
\newtheorem{notation}{Notation}
\newtheorem{proposition}{Proposition} 

 \newcommand{\cqfd}{\hfill
$\square$}

\newcommand{\NN}{\ensuremath{{\bf
N}}}

\newcommand{\dsigma}{\ensuremath{d_{\rho^{\prime}}^{\sigma}}}

\newcommand{\liek}{\ensuremath{\mathfrak{k}}}
\newcommand{\liep}{\ensuremath{\mathfrak{p}}}

\newcommand{\lieg}{\ensuremath{\mathfrak{g}}}
\newcommand{\lien}{\ensuremath{\mathfrak{n}}}
\newcommand{\liel}{\ensuremath{\mathfrak{l}}}
\newcommand{\liea}{\ensuremath{\mathfrak{a}}}
\newcommand{\liem}{\ensuremath{\mathfrak{m}}}
\newcommand{\liez}{\ensuremath{\mathfrak{z}}}
\newcommand{\lieh}{\ensuremath{\mathfrak{h}}}

\newcommand{\X}{\ensuremath{{\bf X}}}

\newcommand{\R}{\ensuremath{{\bf R}}}
\newcommand{\RR}{\ensuremath{{\Bbb R}}}
\newcommand{\CC}{\ensuremath{{\Bbb C}}}
\newcommand{\hh}{\ensuremath{{\Bbb H}}}
\newcommand{\OO}{\ensuremath{{\Bbb O}}}
\newcommand{\KK}{\ensuremath{{\Bbb K}}}

\newcommand{\cB}{\ensuremath{B}}

\newcommand{\cD}{\ensuremath{{\cal D}}}

\newcommand{\hx}{\ensuremath{{\hat x}}}
\newcommand{\hy}{\ensuremath{{\hat y}}}

\newenvironment{preuve}{\begin{trivlist}\item {\it Proof}:}  {\cqfd
\end{trivlist}} \newenvironment{thm*}{\begin{trivlist}\item {\bf Th\'eor\`eme }
\em}{\end{trivlist}} \newenvironment{prop*}{\begin{trivlist}\item {\bf
Proposition } \em} {\end{trivlist}} \newenvironment{coro*}{\begin{trivlist}\item
{\bf Corollaire } \em} {\end{trivlist}}

\newenvironment{lem*}{\begin{trivlist}\item {\bf Lemme } \em}{\end{trivlist}}

\begin{document}
\title{A Ferrand-Obata theorem  for rank one parabolic geometries\\}
\author{Charles Frances}
\date{}
\maketitle
\noindent{{\bf Abstract.}}
 The aim of this article is the proof of  the following result:\\
\ \\
\noindent {\bf Theorem}. 
{\it Let  $M$ be  a connected manifold endowed with a  regular Cartan geometry  $(M,B,\omega)$  modelled on the boundary $\X = \partial {\bf H}_{\KK}^d$ of a $d$-dimensional hyperbolic space  ${\bf H}_{\KK}^d$ over the field $\KK$ ( $\KK= \RR, \CC, \hh$ or the algebra of octonions  $\OO$). If the group of automorphisms $Aut(M, \omega)$ does not act properly on $M$, then $M$ is geometrically isomorphic to:

$\bullet$  $\X$ if $M$ is compact.

$\bullet$  $\X$ minus a point in the other cases.}


\section{Introduction}

At the very begining of the seventies, several works inspired by the so called {\it Lichnerowicz's conjecture} about the conformal group of Riemannian manifolds, led to the following result: 

\begin{theoreme}[Ferrand-Obata]{\cite{ferrand1},\cite{ferrand2},\cite{obata}}
\label{ferrand}

Let  $(M,g)$ be a Riemannian manifold of dimension $n \geq 2$. If the group of conformal transformations of  $(M,g)$ does not act properly on $M$, then $(M,g)$ is conformally equivalent to:

$\bullet$ the standard conformal sphere if $M$ is compact.

$\bullet$ the Euclidean space if $M$ is noncompact.

\end{theoreme}
Let us recall here that the action of a group $G$, by homeomorphisms on a manifold $M$, is said to be   {\it proper} if for every compact subset  $K \subset M$, the set:

\[ G_K=\{  g \in G \ | \ g(K) \cap K \not = \emptyset   \}      \]

has compact closure in $Homeo(M)$ (where $Homeo(M)$, the group of homeomorphisms of $M$ is endowed with the compact-open topology). We do not suppose  {\it a priori} that $G$ is closed in $Homeo(M)$.

Let us precise that the correct and definitive proof in the noncompact case was achieved by  J.Ferrand alone in \cite{ferrand2}. The proof of M.Obata dealt only with the compact case, under the more restrictive asumption that the action of the {\it identity component} of the conformal group, is nonproper. A different proof of Theorem  \ref{ferrand} has been proposed by  R.Schoen in  \cite{schoen}. The methods involved allowed moreover to get a similar result for $CR$ structures:

\begin{theoreme}[Schoen]
\label{schoen}
Let $M^{2n+1}$, $n \geq 1$, be a manifold endowed with a strictly pseudo-convex $CR$ structure. Then , if the group of  $CR$ automorphisms of  $M$ does not act properly on $M$, then:

$\bullet$ $M$ is $CR$-diffeomorphic  to the standard $CR$ sphere if $M$ is compact.

$\bullet$ $M$ is $CR$-diffeomorphic  to the  Heisenberg group, endowed with its standard structure, if $M$ is noncompact.

\end{theoreme}

\subsection{ Cartan geometries}

Although the methods used by J.Ferrand and R.Schoen are completely different, the statements of Theorems \ref{ferrand} and \ref{schoen} are very analogous, and let suspect that they are two aspects of a more general result.  From the geometric point of view, there is a concept which unify, in dimension $n \geq 3$, Riemannian conformal structures and strictly pseudo-convex $CR$ structures:  that of {\it Cartan geometry}. 

Intuitively, a Cartan geometry is the data of a manifold infinitesimally modelled on some homogeneous space $\X=G/P$, where $G$ is a Lie group and $P$ a closed subgroup of $G$. More technically, a Cartan geometry on a manifold $M$, modelled on the homogeneous space $\X=G/P$, is the data of:

\noindent $(i)$ a principal $P$-bundle $B \rightarrow M$ over $M$.

\noindent $(ii)$ a $1$-form $\omega$ on $B$, with values in the Lie algebra $\lieg$, called {\it Cartan connexion}, and satisfying the following conditions:

- At every point $p \in B$, $\omega_p$ is an isomorphism between $T_pB$ and $\lieg$.

- If $X^{\dagger}$ is a vector field of $B$, comming from the action by right multiplication of some one-parameter subgroup $t \mapsto Exp_G(tX)$ of $P$, then $\omega(X^{\dagger})=X$.

- For every $a \in P$, ${R_a}^*\omega=Ad(a^{-1})\omega$ ($R_a$  standing for the right action of $a$ on $B$).

A Cartan geometry on a manifold $M$ will be denoted by  $(M,B,\omega)$.

A lot of classical geometric structures can be interpreted in terms of Cartan geometry. The most famous ones are Riemannian metrics  (resp. pseudo-Riemannian metrics of signature $(p,q)$). In this case, the space $\X$ is just the Euclidean space (resp. the Minkowski space of signature $(p,q)$), $G$ is the group of isometries $SO(n)\ltimes \RR^n$ (resp.  $SO(p,q) \ltimes \RR^n$), and $P$  the normal subgroup constituted by the translations.

Since the definition of a Cartan geometry is not very intuitive, two natural problems arise at once. The first is the interpretation of the data of a Cartan connection $\omega$ on a principal bundle $B \rightarrow M$, as an underlying geometric structure on $M$. In a lot of interesting geometric situations, such an interpretation is available (see for instance \cite{cap1}, \cite{morimoto}, \cite{tanaka2}, and references therein). The second interesting problem is to know, if a given underlying structure on a manifold $M$, determines canonically a Cartan geometry. This problem, known as the {\it equivalence problem}, is  generally quite difficult. It was solved by E.Cartan himself for conformal Riemannian structures, in dimension $n \geq 3$, and for strictly pseudo-convex $CR$ structures  (\cite{cartan1} in dimension $3$, \cite{tanaka1} and \cite{chern} for dimensions $n \geq 3$. See also \cite{kobayashi} and \cite{sharpe} for the conformal case).  Otherwise stated, if a manifold $M$ of dimension $n \geq 3$, is endowed with a conformal class of Riemannian metrics (resp. with a strictly pseudo-convex $CR$ structure), one is able to build a $P$-principal bundle $B$ over $M$, and a Cartan connection $\omega$ on it. In the conformal (resp. $CR$) case, the model space $\X$ is the conformal sphere ${\bf S}^n= \partial {\bf H}_{\RR}^{n+1}$ (resp. the $CR$ sphere ${\bf S}^{2n+1}=\partial {\bf H}_{\CC}^{n+1}$), the group $G$ is the Moebius group $SO(1,n+1)$ (resp. the group $SU(1,n+1)$), and $P$ is a parabolic subgroup: the stabilizer of a point on $\X$. Moreover, if one requires that it satisfies suitable normalizations conditions, the Cartan connection $\omega$ is defined uniquely. Thus, any conformal diffeomorphism (resp. $CR$ diffeomorphism) acts on $B$ preserving the connection $\omega$.

For any Cartan geometry $(M,B,\omega)$, we define $Aut(B, \omega)$ as the set of $C^1$-diffeomorphisms $\phi$ of $B$, such that $\phi^* \omega=\omega$. Every $\phi$ of $Aut(B, \omega)$ commutes with the right action of $P$ on $B$, so that $\phi$ induces a diffeomorphism $\overline \phi$ of $M$. The subset of diffeomorphisms of $M$ obtained in this way is denoted by $Aut(M, \omega)$.

There is also a notion of geometric equivalence for two Cartan geometries $(M, B, \omega)$ and $(N, B^{\prime}, \omega^{\prime})$, modelled on the same space $\X=G/P$. We say that $M$ and $N$ are {\it geometrically isomorphic}, if there is a diffeomorphism  $\phi : B \rightarrow B^{\prime}$ such that $\phi^* \omega^{\prime}= \omega$.

\subsection{Statement of results}

The aim of the  article is the generalization of Theorems \ref{ferrand} and \ref{schoen} to any Cartan geometry modelled on spaces $\X=G/P$, where $G$ is a simple Lie group of real rank one, with finite center, and $P$ is a parabolic subgroup of $G$. These spaces $X$ are the boudaries of the different hyperbolic spaces  ${\bf H}_{\KK}^d$,  $\KK$ standing for  the field $\RR$ of real numbers, $\CC$ of complex numbers, $\hh$ of quaternions,  or the octonions $\OO$. We will make the asumption  $d \geq 2$ if $\KK=\RR$, and $d \geq 1$ otherwise (except for the octonionic case, where the dimension $d$ is necessarily 2).  
Implicitely, when we will speak about a Cartan geometry modelled on  $\partial {\bf H}_{\KK}^d$,  we will always see   $\partial {\bf H}_{\KK}^d$ as the homogeneous space $G/P$, where $G=Iso({\bf H}_{\KK}^d)$, and $P$ is the stabilizer, in  $Iso({\bf H}_{\KK}^d)$, of a point of  $\partial {\bf H}_{\KK}^d$.

Let us recall the groups $G$ involved:

- $G=SO(1,n)$, $d \geq 2$ for $\KK=\RR$. The space  $\X=\partial {\bf H}_{\RR}^d$ is a sphere ${\bf S}^{d-1}$.

- $G=SU(1,d)$, $d \geq 1$ for $\KK=\CC$.  The space $\X=\partial {\bf H}_{\CC}^d$ is a  sphere ${\bf S}^{2d-1}$.

- $G=Sp(1,d)$, $d \geq 1$ for $\KK=\hh$.  The space $\X=\partial {\bf H}_{\hh}^d$ is a sphere ${\bf S}^{4d-1}$.

- $G=F_4^{-20}$ if $\KK=\OO$.  The space $\X=\partial {\bf H}_{\OO}^2$ is a sphere ${\bf S}^{15}$.\\
\ \\

We can now state our main result:

\begin{theoreme}
\label{theoreme1}

Let $(M,B,\omega)$ be a Cartan geometry modelled on $\X=\partial {\bf H}_{\KK}^d $, the boundary of the $d$-dimensional hyperbolic space over  $\KK=\RR, \CC, \hh$ or $\OO$. We suppose that $M$ is connected, and that the connexion $\omega$ is regular. Then if  $Aut(M,\omega)$ does not act properly on $M$, then $M$ is geometrically isomorphic to:

$\bullet$ $\X$ if $M$ is compact.

$\bullet$  $\X$ minus a  point otherwise.
\end{theoreme}

The hypothesis of regularity is a technical hypothesis on the curvature of the connection $\omega$ (which is mild, since satisfied in most of interesting cases). The notion of regularity will be explained in section \ref{cartan}.

The conclusions of Theorem \ref{theoreme1} will hold for any geometric structure on a manifold $M$, for which the equivalence problem has been solved, and from  which one is able to define a {\it regular} canonical Cartan geometry, modelled one one of the spaces $\partial {\bf H}_{\KK}^d$. Since the geometries involved in Theorem \ref{theoreme1} are parabolic geometries, we can use the  works done on the equivalence problem for these geometries. We get that the theorem will apply to  conformal Riemannian structures, and strictly pseudo-convex $CR$-structure in dimension $n \geq 3$ (this gives a unified proof for Theorems \ref{ferrand} and \ref{schoen}), but also to partially integrable almost $CR$-structures (\cite{cap1}, \cite{morimoto}, \cite{tanaka2}), as well as to quaternionic and octonionic contact structures introduced by O.Biquard (\cite{biquard}, see also \cite{cap3}). 

\subsection{Ideas of the proof and organisation of the paper}

The proof of Theorem \ref{theoreme1} is based essentially on the understanding of the dynamics of sequences of automorphisms of a manifold $M$, endowed with a Cartan geometry modelled on some $\partial {\bf H}_{\KK}^d $.   The main point is to prove that if a sequence $(f_k)$ of $Aut(M, \omega)$ does not act properly on $M$, it has the property:

\noindent{$(P)$: {\it there is an open subset $U \subset M$ which collapses  to a point under the action of $(f_k)$. }}

This is a fundamental property since one can prove that it implies the flatness of the geometry on the open subset $U$ (see Proposition \ref{courbure}).

In fact, J.Ferrand and R.Schoen proved also the property $(P)$ for the sequences of conformal (resp. $CR$) diffeomorphisms which do not act properly. Let us observe that they both did it using analytical tools (J.Ferrand writes in the introduction of  \cite{ferrand2} "{\it In fact,} [Theorem \ref{ferrand}] {\it is not actually concerning the theory of Lie groups and may be considered as a mere theorem of Analysis}").

The way we adopt to prove the property $(P)$ is on the contrary purely geometric. Let us begin, recalling that this property $(P)$ is typical of the sequences of $Iso({\bf H}_{\KK}^d)$, when they act on the boundary $\X=\partial {\bf H}_{\KK}^d$. That is what is usually called a dynamics of {\it convergence type} (we say also  "North-South" dynamics;  see section \ref{dynamique}). The main idea of Theorem \ref{theoreme1} is to use the Cartan connection to link the dynamical  properties  of sequences of $Aut(M, \omega)$, to dynamical properties  of sequences of $G$, acting on $\X$. Let us begin with the simple case of a flat Cartan geometry. Such a geometry is more comonly called  a $(G, \X)$-structure on $M$. In this case, one can define a developping map $\delta: \tilde M \to \X$, which is an immersion, as well as a morphism $\rho : Aut(\tilde M, \tilde \omega) \to G$ (we refer for example to \cite{thurston} for general results on $(G,\X)$-structures). Moreover, the equivariance relation $\delta \circ f=\rho(f) \circ \delta$ is satisfied for every $f \in Aut(\tilde M, \tilde \omega)$. This equivariance relation is crucial since it allows to recover, at least locally, the dynamics of a sequence $(f_k)$ of $Aut(\tilde M, \tilde \omega)$ from the dynamics of $\rho(f_k)$ on $\X$ (see, for example,  \cite{lafontaine} or \cite{frances-tarquini}  as an illustration).

Of course, the Cartan geometry we are looking at is generally not flat, {\it a priori}. Although all the tools used in the case of $(G, \X)$-structures break down in this case, it still remains something of the previous scheme.  Let us fix $x_0 \in M$, $\hx_0 \in B$ over $x_0$, and $o \in \X$. The Cartan connection still allows to define some kind of developping map, denoted ${\cal D}_{x_0}^{\hx_0}$. This is a map from the space of curves of $M$ passing through $x_0$, to the space of curves of $\X$ passing through $o$. This (classical) procedure will be recalled in section \ref{developpante}.

Now, let $(f_k)$ be a sequence of $Aut(M, \omega)$. To simplify and avoid technicalities, we suppose here that $(f_k)$ fixes $x_0$. In section \ref{holonomie}, we will explain how  to associate to $(f_k)$ some {\it holonomy sequence} $(b_k)$ of $P$. 
The fundamental point is that  we still have some equivariance relation:

\begin{equation}
\label{equi}
{\cal D}_{x_0}^{\hx_0} \circ f_k = b_k \circ {\cal D}_{x_0}^{\hx_0}
\end{equation}

Let us insist on the fact that this is a general principle, which will be probably useful for the dynamical study of automorphisms of other Cartan geometries than those of this paper. Relation (\ref{equi})  shows a link between the action of $(f_k)$ on the curves of $M$ passing through $x_0$, and the action of $(b_k)$ on the curves of $\X$ passing through $o$. Of course, the space of curves of $\X$ (resp. of $M$) passing through $o$ (resp. through $x_0$) is to much complicated, and we will restrict ourself to 
 the action on a class of curves, which are distinguished for the geometry under consideration: the geodesics of the Cartan geometry. These geodesics are defined in section \ref{geodesique}. They coincide with the conformal geodesics when $\X= \partial {\bf H}_{\RR}^d$, and the chains introduced by E.Cartan when $\X =\partial {\bf H}_{\CC}^d $.
 
 The rough idea to get the property $(P)$, is that there will still be   some trace of the  "North-South" dynamics, for the action of $(b_k)$ on the set of geodesics passing through $o$ (that is the result of the study done in section \ref{dynamique segments}). Thanks to the relation (\ref{equi}), this dynamical behaviour will be transmitted to the action of $(f_k)$ on the set of geodesics passing through $x_0$. We then recover the dynamics of $(f_k)$ around $x_0$ by some kind of  projection.  
 
 Let us remark that the "North-South" dynamical behaviour caracterizes the rank one situation. This is basically why Theorem \ref{theoreme1} has no analogous for certain  higher rank parabolic geometries (see for example \cite{frances2}, which deals with the conformal Lorentzian situation).
 
 As already said, once the property $(P)$ is proved, we get the flatness of the Cartan geometry on an open subset $U \subset M$. To get Theorem \ref{theoreme1}, it still remains to show that $U$ is in fact the whole $M$, or $M$ minus a point. This is done thanks to a rigidity result for geometrical embeddings of certain Cartan geometries, of independant interest.  This result is stated at the end of section \ref{demonstrations} (Theorem \ref{plongement}), and the last part of the article is devoted to its proof.


\section{Geometry of the model spaces  ${\bf X}$}
\label{geox}
\subsection{Algebraic preliminaries}
\label{algebriques}

Most of the following preliminaries are very clearly exposed in \cite{knapp}.

Let $G$ be a simple Lie group of real rank $1$, with finite center , and  $\lieg$ its Lie algebra.     We choose a Cartan involution $\theta$   on $\lieg$. This involution determines a Cartan decomposition $\lieg = \liek_{\theta} \oplus \liep_{\theta}$.     Let $\liea $ be a maximal abelian subalgebra of   $\liep_{\theta}$. Since we supposed the rank of $G$ to be $1$, the dimension of $\liea$ is also $1$.  Let $\Delta$ be the set of roots for the adjoint representation of $\liea$ on $\lieg$.  For any  $\lambda \in \Delta$, the space  $\lieg_{\lambda}$ is defined as  $\lieg_{\lambda} = \{v \in \lieg \ | \ Ad(a^t)v=e^{t\lambda(X_0)}v   \}$. Since the rank of $G$ is one, there are two possibilities for $\Delta$ :

- $\Delta = \{ -\alpha, + \alpha  \}$ (case $\lieg=\mathfrak{so}(1,n)$, $n \geq 2$, or $\mathfrak{su}(1,1)$ or $\mathfrak{sp}(1,1)$).

- $\Delta = \{ -2 \alpha, -\alpha, + \alpha , +2 \alpha \}$ (all the other cases). 

We choose      $X_0 \not = 0$ in  $\liea$ such that $\alpha(X_0)=1$, and denote by $A$   the one parameter subgroup $a^t=Exp_G(tX_0)$.  


In any of the two cases, the Lie algebra $\lieg $ admits the decomposition  $\lieg=\lien^- \oplus \liea \oplus \liem \oplus \lien^+$. The algebra $\liel=\liea \oplus \liem$ is the centralizer of $\liea$ in $\lieg$. It is stable under the action of the involution $\theta$, and $\liea$ is the eigenspace of $\liel$ associated to the eigenvalue $-1$, since $\liem$ is just the space of fixed points of $\theta$ on $\liel$.   

One can write $\lien^+=\lien_1^+ \oplus \liez^+$ (resp.  $\lien^-=\liez^- \oplus \lien_1^-$) where $\liez^+$ (resp. $\liez^-$) is the center of  $\lien^+$ (resp. $\lien^-$). 

 Precisely, if  $\Delta = \{  - \alpha, + \alpha\}$, we simply have $\liez^+=\lien^+=\lieg_{+ \alpha}$ (resp. $\liez^-=\lien^-=\lieg_{- \alpha}$), and $\lien_1^+=\lien_1^-=\{  0\}$.

When $\Delta = \{ -2 \alpha , - \alpha, + \alpha, +2 \alpha \}$, then $\lien^+=\lieg_{+ \alpha} \oplus \lieg_{+ 2 \alpha}$ (resp. $\lien^-=\lieg_{-2 \alpha} \oplus \lieg_{- \alpha}$), $\liez^+=\lieg_{+ 2 \alpha}$ (resp.  $\liez^-=\lieg_{- 2 \alpha})$) and $\lien_1^+=\lieg_{+ \alpha}$ (resp. $\lien_1^-=\lieg_{- \alpha}$).

The two Lie algebras $\lien^+$ et $\lien^-$ are nilpotent. The exponential $Exp_G$ is a diffeomorphism from  $\lien^+$ (resp. $\lien^-$) onto a closed subgroup   $N^+ \subset G$ (resp.  $N^- \subset G$).

The parabolic subgroup $P$ is the closed subgroup of  $G$ with Lie algebra $\liep = \liea \oplus \liem \oplus \lien^+$. We denote  $\pi_{X}$ the  projection $G \to G/P= \X$. As we already said in the introduction, the space $\X$ is diffeomorphic to a sphere.

Let us also recall the  {\it Langland's decomposition} $P=MAN^+$ for the group $P$, where $M$ and $A$ are closed subgroups of $G$, with respective Lie algebras $\liem$ et $\liea$. We will denote $L=MA$.

\subsection{Charts}
\label{bruhat}
{- \it Atlas :}
 Let us call $o$ the projection of   $P$ on  $G/P$.  As the rank of $G$ is one, its Weyl group is reduced to two elements: $W(G)= \{  e,w \}$. The element $w$ acts on  $\X$ by an  involution, and sends $o$ to a point  $\nu \in \X$. The  points $o$ and  $\nu$ are left fixed by the group $L$. Also, $wN^+w=N^-$. The Bruhat decomposition (see \cite{knapp} for example) writes $G=P \cup wN^+wP=P \cup N^-P$. In other words, if one calls $\Omega_{o}$ (resp. $\Omega_{\nu}$) the orbit of $o$ under the action of  $N^-$ (resp. the orbit of $\nu$ under the action of  $N^+$), the manifold $\X$ can be written as the union : $\X=\{  o \} \cup \Omega_{\nu}= \{  \nu \} \cup \Omega_{o}$.

This gives an atlas with two charts on $\X$ : the mapping  $s^- : \lien^- \to {\bf X}$ defined by $s^-(u)=Exp_G(u).o$, which is a diffeomorphism from  $\lien^-$ onto the open set  $\Omega_{o}$. We also have the chart  $s^+ : \lien^+ \to {\bf X}$ defined by $s^+(u)=Exp_G(u).\nu$, which maps $\lien^+$ diffeomorphically onto $\Omega_{\nu} $. \\

\noindent{- \it Auxiliary metrics:}
 we endow $\lieg$ with a scalar product  $< \ >_{\lieg}$, invariant by the Cartan involution $\theta$. We denote by  $||.||$ the associated norm on $\lieg$. We extend this scalar product to a left invariant Riemannian metric $\rho_G$ on  $G$. This metric induces on  $N^+$ and  $N^-$ two left invariant Riemannian metrics  $\rho^+$ and $\rho^-$, that we carry on  $\Omega_{\nu}$ and  $\Omega_{o}$. We thus get two Riemannian metrics on $\Omega_{\nu}$ and  $\Omega_{o}$ (that we still denote by $\rho^+$ and $\rho^-$), for which the actions of $N^+$ and  $N^-$ are isometric.\\


\noindent{- \it Change of charts} :
We denote  $s_-^+=(s^+)^{-1} \circ s^-$. The map $s_-^+$ maps $\lien^- \setminus \{  0\}$ on $\lien^+\setminus \{  0\}$. We are going to give a formula for the restriction of $s_-^+$ to $\liez^- \setminus \{  0\}$.

Let  $u \in \liez^-$, $u \not = 0$, and  $w=[u,\theta u]$. Since we saw that there is  $ \lambda \in \Delta$ such that  $\liez^-=\lieg_{\lambda}$, we get that  $\RR.u \oplus \RR. w \oplus \RR. \theta u$ is a subalgebra of $\lieg$, isomorphic to  ${\mathfrak sl}(2,\RR)$. In fact, we can choose a normalization $u^{\prime}=\lambda u$, $w^{\prime}= \mu w$, and $v^{\prime}=\theta u$, such that $[w^{\prime},u^{\prime}]=2u^{\prime}$, $[w^{\prime},v^{\prime}]=-2 v^{\prime}$ and  $[u^{\prime},v^{\prime}]=w^{\prime}$. The formula of the stereographic projection in dimension one yields a real $a_u$ (which  depends only of the direction of  $u$) such that : 
\begin{equation}
\label{changement}
s_-^+  (u) =   a_u .\frac{\theta u}{||u||^2}
\end{equation}

\subsection{Geodesics on ${\bf X}$}

We begin with some notations. For $u \in \lieg$, we call $u^*$ the curve from $[0,1]$ to ${\X}$ defined by $u^*(t)=\pi_X \circ Exp_G(tu)$. By $[u]$, we will mean the geometrical arc supporting   $u^*$.

\subsubsection{Parametrized geodesics}

We define the following subset  ${\bf Q} \subset \lieg$:

\[   {\bf Q} = \{  v \in \lieg \   |  \  v=Ad(b).u, \ u \in \liez^-, \ b \in P  \}\]

\begin{definition}[Parametrized geodesics and geodesic segments]

We call geodesic of  $\X$ any curve from $I$ to $\X$, where $I$ is an interval of $\RR$ containing $0$, which is of the form  $t \to g.\pi_X \circ Exp_G(tu)$ with  $u \in {\bf Q}$ and $g \in G$. When $I=\RR$, we speak of {\it maximal geodesic}.

- For $u \in {\bf Q}$, the curve $u^*$ is called the parametrized geodesic segment with origin   $o$ associated to $u$.

-  One says that  $[u]$ is  the (unparametrized) geodesic segment with origin $o$, associated to $u$.  The space of geodesic segments with origin $o$ is denoted by ${\bf {[Q]}}$. 

- We will denote by  ${\bf \dot Q}$  (resp. $[\bf \dot Q ]$) the space  ${\bf Q}$ with $0_{\lieg}$ removed (resp. the space $\bf [Q] $ with the "trivial" segment  $[o]$ removed).
\end{definition}

In all what follows, we will endow ${\bf [Q]}$ with the Hausdorff topology on closed sets of $\X$.

\subsubsection{Fundamental properties of ${\bf [Q]}$}

Our first task is to understand the action of $P$ on ${\bf Q}$, and for that, we have to desribe ${\bf [Q]}$ more precisely.

\begin{lemme}
There is a morphism  $\rho : P \to Aff(\lien^+)$ such that for all $b \in P$ : $b \circ s^+=s^+ \circ \rho(b)$.
\end{lemme}

One looks separately at the action of the different components of  Langland's decomposition.

$\bullet$ {\it Action of $L$.} 

Let  $l$ be an element of $ L$, and  $u \in \lien^+$. We set $n^+=Exp_G(u)$ and  $x=n^+.\nu$. Then $l.x=Ad(l)(n^+).l.\nu$, and since   $l. \nu=\nu$, $l.x=s^+(Ad(l).u)$. So, $\rho(l)=Ad(l)_{|\lien^+}$ : the action of  $L$ in the chart $s^+$ is just the adjoint action.\\

$\bullet$ {\it Action of $N^+$.}  

Let $n_0^+=Exp(u_0)$ and  $n^+=Exp(u)$ be two elements of $N^+$. Since $\lien^+$ is a one, or two-step nilpotent Lie algebra, the Campbell-Hausdorff  formula yields $n_0^+n^+=Exp(u_0+u+\frac{1}{2}ad(u_0)(u))$. We thus get that  $n_0^+.s^+(u)=s^+((Id +\frac{1}{2}ad(u_0)).u+u_0)$, or :
\[ \rho(n_0^+) : u \mapsto (Id + \frac{1}{2}ad(u_0)).u + u_0 \]

\begin{remarque}
\label{translations}
If $u \in \liez^+$, $\rho(n_0^+).u=u+u_0$, so that $\rho(N^+)$ acts by  translations on  $\liez^+$.

\end{remarque}

For every half-line $\alpha$ of $\lien^+$, there is a unique point $x \in \lien^+$ and a  unique vector $v \in \lien^+$ of norm $1$ for $||. ||$, such that  $\alpha = \{ y\in \lien^+ \ | \ y=x+tv, \ , t \in \RR_+   \}$. So, if $S_{\lien^+}$ (resp.  $S_{\liez^+}$) is the unit sphere of $\lien^+$ (resp.  $\liez^+$) for the  norm $||.||$, the space of half-lines of $\lien^+$ (resp. of half-lines whose direction is in  $\liez^+$) is identified with the product  $\lien^+ \times S_{\lien^+}$ (resp. $\lien^+ \times S_{\liez^+}$). In the following, such an half-line will be denoted by  $[x,u)$, with $x\in \lien^+$ and  $u \in  S_{\lien^+} $ (resp. $u \in  S_{\liez^+} $).

\begin{remarque}
\label{direction}
The group  $\rho(P)$ leaves globally invariant the set of half-lines of  $\lien^+$ whose direction is in  $\liez^+$.
\end{remarque}
\ \\
\begin{proposition}
\label{fermeture}
\ \\
\noindent $(i)$ Under the map $\mu : [u] \to  {(s^+)}^{-1}( [u] \cap \Omega_{\nu})$, the space ${\bf [\dot Q]}$ is homeomorphic to the space of affine half-lines of $\lien^+$  whose direction is in $\liez^+$.

\noindent $(ii)$There exists a continuous section  $s : {\bf [ \dot Q]} \to {\bf \dot Q} $.
\end{proposition}

\begin{preuve}

Let us first remark that for all $b \in P$, the following equivariance relation is true :

\[ b.[u]=\rho(b).\mu([u])  \]

Then, the formula (\ref{changement})  ensures that if  $u \in \liez^-$  then  $\mu([u])$ is the half-line  $[a_u \frac{\theta u}{||u||^2},  \frac{\theta u}{||u||})$. Now, every  $u^{\prime} \in {\bf Q}$ writes  $Ad(b).u$ for $b \in P$ and  $u \in \liez^-$. From the relation $ b.[u]=\rho(b).\mu([u])$ and the remark  \ref{direction}, one infers that the image of $\mu$ is included in the set of affine half-lines of  $\lien^+$, whose direction is in  $\liez^+$. Since  $\mu$ is clearly an homeomorphism on its image, we only have to show that $\mu$ is surjective to conclude the proof. From the formula (\ref{changement}), it is clear that any half-line $[u, u)$, with  $u \in {S}_{\liez^+}$,  is in the image of $\mu$. Then, we make $\rho(N^+)$ act on these half-line, and we get that all the half-lines of   $\lien^+$, whose direction is in  $\liez^+$ are in the image of $\mu$.

\ \\
We still have to define the section $s$. We define first $\tilde s : \lien^+ \times S_{\liez^+} \to {\bf \dot Q}$ by : $\tilde s : [x,v) \mapsto Ad(Exp_G(x-v)).\frac{\theta v}{a_{\theta v}}  $

Then, we set  $s=\tilde s \circ \mu$. Let us check that  $s$ is really a section:

$[\tilde s([x,v))]=Exp_G(x-v).[\frac{\theta v}{a_{\theta v}}]
=\mu^{-1}(\rho(Exp_G(x-v)).[v,v))=\mu^{-1}([x,v))$.

\end{preuve}

\begin{lemme} 
\label{pro}

The space  $\bf [Q]$ has the following properties :

\noindent $(P_1)$  For every point $x \in \X$, there is a segment $\alpha \in \bf [Q]$  joining $o$ and $x$.

\noindent $(P_2)$  For the  Hausdorff topology, $\bf [Q]$ is a closed subset of ${{\cal K}(\Omega_o)}$, the set of compact subsets of $\Omega_o$.

\noindent $(P_3)$ Let $(\alpha_k)$ be a sequence of $\bf [Q]$ tending to $[o]$, then $\lim_{k \to + \infty}L^-(\alpha_k)=0$, where $L^-(\alpha_k)$ is the length of the segment  $\alpha_k$ for the metric $\rho^-$.

 \end{lemme}

\begin{preuve}

Property $(P_1)$ is trivial if  $x=o$. If it is not the case, $x \in \Omega_{\nu}$. We set $y=(s^+)^{-1}(x)$ and we choose $u \in S_{\liez^+}$. Then $\mu^{-1}([y,u))$ is a segment of  ${\bf [Q]}$ joining  $o$ and $x$.

Property $(P_2)$ follows from the homeomorphism   $\mu$ between  ${\bf [{\dot Q}]}$ and the set of half-lines of  $\lien^+$ whose direction is in  $\liez^+$, and the fact that, up to subsequence, any diverging sequence of $\lien^+ \times S_{\liez^+}$ either tends to $[o]$, or leaves every compact subset of ${\cal K}(\Omega_o)$. 

  We now prove property  $(P_3)$. Let us suppose tat $(P_3)$ is not true, and suppose the existence of a sequence $\alpha_k \in [\bf {\dot Q }]$ which tends to  $[o]$, and such that $L^-(\alpha_k)> \epsilon>0$, for all $k \in \NN$. We identify  $[\bf {\dot Q }]$ with the space of half-lines in $\lien^+$ whose  direction is in  $\liez^+$, so that we note $\alpha_k=[x_k,u_k)$, with $x_k \in \lien^+$ and $u_k \in S_{\liez^+}$. Looking if necessary at a subsequence, we will suppose that $u_k$ has a limit  $u_{\infty} \in S_{\liez^+}$.  In the proof, we will use the following criteria, whose proof is easy:

\begin{lemme}
\label{outil}
 Let  $[x_k,u_k)$ be a sequence of  ${\bf [{\dot Q}]}$. We suppose that  $x_k \to \infty$ (i.e $(x_k)$ leaves every compact subset of  $\lien^+$), and that there are $u_{\infty}$ and $v_{\infty}$ in  $S_{{\liez}^+}$ such that $\lim_{k \to +\infty}u_k=u_{\infty}$ and  $\lim_{k \to +\infty}\frac{x_k}{||x_k||}=v_{\infty}$. Then the sequence $[x_k,u_k)$ tends to $[o]$ if and only if  $v_{\infty} \not = - u_{\infty}$.

\end{lemme}

Let $a^t$  be the  Cartan flow introduced in section \ref{algebriques}. The flow $a^t$ acts on $\lien^-$ by transformations of the form $$\left( \begin{array}{cc} e^{-2t}Id_{\liez^-} & 0 \\0 & e^{-t}Id_{\lien_1^-}  \end{array}  \right) $$  So, the norm associated to any scalar product on  $\lien^-$ is contracted exponentially by $ad(a^{t})$ (when $t \geq 0$). As a consequence, $a^{t}$ acts also by contractions for the metric $\rho^-$ : there is a constant  $c >0$ such that  $L^-(a^{t}.[\sigma]) \leq e^{-ct} L^-([\sigma])$.

For any sequence  $(x_k)$ tending to infinity in  $\lien^+$ there is a sequence $(t_k)$ of reals, with $\lim_{k \to + \infty}t_k = - \infty$ such that  $Ad(a^{t_k}).x_k$ remains in a compact subset of $\lien^+ \setminus \{  0_{\lien^+} \}$. Looking if necessary at a subsequence, we will suppose that $\sigma_k=a^{t_k}.\alpha_k$ tends to  $\sigma_{\infty}=[x_{\infty},u_{\infty})$, with $x_{\infty} \not = 0_{\lien^+}$.   

If there is no $\mu>0$ such that  $x_{\infty}  = - \mu u_{\infty}$, the segment $\sigma_{\infty}$ does not contain $0_{\lien^+}$, and thus  $L^-([\sigma_k])$ is bounded from above by a real  $M$, for any $k \in \NN$.    This yields a  majoration of the form $L^-(\alpha_k) \leq Me^{ct_k}$, which contradicts the hypothesis  $L^-(\alpha_k)>\epsilon$ (because $t_k \to - \infty$).

Now, if there exists  $\mu>0$ such that  $x_{\infty}= - \mu u_{\infty}$. We write $x_k=y_k+z_k$, with  $y_k \in \lien_1^+$ and  $z_k \in \liez^+$. One has  $\lim_{k \to +\infty} e^{t_k}y_k=0$ and  $\lim_{k \to + \infty}e^{2 t_k}z_k=x_{\infty}=-\mu u_{\infty}$. So, $\lim_{k \to + \infty}e^{2 t_k}(y_k+z_k)=- \mu u_{\infty}$, which implies  $\lim_{k \to + \infty} \frac{x_k}{||x_k||}=- u_{\infty}$. But lemma \ref{outil} then ensures that the sequence  $[x_k,u_k)$ did not tend to $[o]$, yielding a  contradiction.

\end{preuve}

\subsection{Dynamical aspects}
\label{dynamique}

\subsection{Some general definitions }
We now introduce some dynamical notions which will be used all along the article. 

Let $M$ be a manifold, and  $G$ a subgroup of homeomorphisms of $M$. 

\begin{definition}[Stable data]

We call stable data  a quadruple $(f_k,x_k, x_{\infty}, y_{\infty})$, where $(f_k)$ is a sequence of $G$,  $(x_k)$ is a sequence of  $M$ converging to $x_{\infty} \in M$, and such that the sequence $y_k=f_k(x_k)$ converges to $y_{\infty} \in M$. 
\end{definition}

It is easy to check that a subgroup $G \subset Homeo(M)$ does not act properly (i.e has a sequence which does not act properly) if and only if there exists a stable data  $(f_k,x_k,x_{\infty},y_{\infty})$, where $(f_k)$ is a sequence of $G$ tending to infinity in  $Homeo(M)$ (we say that a sequence tends to infty in  $Homeo(M)$ if its intersection with any compact subset of $Homeo(M)$ is finite).

\begin{definition}[Equicontinuity]

Let $(f_k,x_k,x_{\infty},y_{\infty})$ be a stable data of $G$. One says that the action of  $(f_k)$ is equicontinuous at  $x_{\infty}$ if there exists a subsequence  $(f_k^{\prime})$ of $(f_k)$, such that for any sequence $x_k^{\prime}$ tending to $x_{\infty}$, $f_k^{\prime}(x_k^{\prime})$ tends to $y_{\infty}$.

\end{definition}

\subsubsection{North-South dynamics on  $\X$}
\label{NS}

\begin{lemme}
\label{nord-sud}
Let $(g_k)$ be a sequence of $G$ tending to infinity. Then, looking if necessary at a subsequence of $(g_k)$, there exist two points $o^+$ and  $o^-$ of ${\bf X}$ (which can be the same), such that $(g_k)$ has the following dynamical properties:

\noindent $(i)$ For every compact subset $K \subset \Omega_{o^+}=\X \setminus \{  o^- \}$, $\lim_{k \to + \infty} g_k(K)=o^+$.

\noindent $(ii)$ For every compact subset $K \subset \Omega_{o^-}=\X \setminus \{  o^+ \}$, $\lim_{k \to + \infty} (g_k)^{-1}(K)=o^-$.
\end{lemme}

\begin{preuve}
The group $G$ admits a  {\it Cartan's decomposition }  $G=KAK$, where $K$ is the maximal compact subgroup of $G$. Thus, it is sufficient to do the proof of the lemma for a sequence  $g_k=a_k$ of elements of  $A$, tending to infinity. Considering if necessary  $a_k$ instead of  $a_k^{-1}$, we get $a_k=a^{t_k}$, with $t_k \to + \infty$.  The matricial expression of the action of  $Ad(a^t)$ on $\lien^-$, given in the proof of  lemma \ref{pro} shows that for any compact subset $K$ of  $\lien^-$, $Ad(a^t).K$ tends uniformly to $0$ as  $t \to + \infty$. Now, from the relation $s^-(Ad(a^{t_k}).u)=a^{t_k}.s^-(u)$, we get the  point $(i)$ of the lemma, setting  $o^+=s^-(0)$. The point $(ii)$ is proved doing the same work on $a^{-t_k}$,  acting on $\lien^+$, and setting  $o^-=s^+(0)$.
\end{preuve}

As a consequence of lemma  \ref{nord-sud} we get the:
\begin{corollaire}
\label{kak}
Let $(b_k)$ be a sequence of  $P$ tending to infinity, with  poles $o^+$ and $o^-$. The action of  $(b_k)$ is equicontinuous at  $o$ if and only if $o=o^+$, and  $o^- \not = o$. In this case, there exist   $l_{1,k}$ and  $l_{2,k}$  two sequences of $P$,  relatively compact in  $P$,  such that $b_k=l_{1,k}a_kl_{2,k}$, where $a_k \in A$.
\end{corollaire}

\begin{preuve}
If $(b_k)$ is a sequence of  $P$ tending to infinity, the dynamical study made in the lemma proves that the only possible fixed points for  $(b_k)$ are $o^+$ and  $o^-$. Also, the action of $(b_k)$ is never equicontinuous at $o^-$ (but it is equicontinuous on the whole  $\Omega_{o^+}$). We infer that $o=o^+$, and  $o^- \not=o^+$. Now, thanks to the previous lemma, $\lim_{k \to + \infty}b_k^{-1}.\nu=o^-$. Since $o^- \in \Omega_{\nu}$, there is a sequence $n_k$ in $N^+$, relatively compact in  $N^+$, so that  $n_kg_k^{-1}.\nu=\nu$. The sequence  $n_kg_k^{-1}$ fixes $o$ and $\nu$, so that it is a sequence of  $L=AM$. We deduce the existence of a sequence $m_k$ of  $M$ such that $n_kg_k^{-1}m_k$ is in $A$, which concludes the proof.
\end{preuve}

\subsubsection{Dynamics of $P$ on the space  ${\bf [Q]}$}
\label{dynamique segments}

\begin{proposition}
\label{dynamique}
Let $(b_k)$be a sequence of  $P$ which tends to infinity. Then, considering if necessary  $(b_k^{-1})$ instead of $(b_k)$, and looking if necessary at a subsequence: 

$(i)$ There is an open set ${\bf \Omega}^+ \subset  {\bf [\dot Q]}$, containing $[o]$ in its closure,  such that for any compact subset $K$ of ${\bf \Omega}^+$, $b_k.K$ tends to $[o]$ as $k$ tends to $+ \infty$.

$(ii)$ The closure of the set  $s(\Omega^+)$ contains elements of  $\liez^- \setminus \{ 0_{\lieg} \}$ which are arbitrarily close to $0_{\lieg}$.

\end{proposition}

\begin{preuve}
The  point $(ii)$ is just a technical point that will be useful in section \ref{demonstrations}.

In the whole proof, we identify  ${\bf [{\dot Q}]}$ with $\lien^+ \times S_{\liez^+}$ thanks to the map $\mu$. 
We don't change the conclusions of Proposition  \ref{dynamique} if we compose $(b_k)$ by a sequence of the compact group $M \subset P$ (observe that $Ad(M)$ preserves $\liez^-$). So, we will suppose that $(b_k)$ is a sequence of $AN^+ \subset P$, and we write $b_k=a^{t_k}n_k^+$. On $\lien^+$, $\rho(b_k)=L_k+T_k$, namely $\rho(b_k)$ is the composition of the linear map $L_k$ and the translation of vector $T_k$. In a basis compatible with the grading $\lien^+=\lien_1^+ \oplus \liez^+$, the matrix of $L_k$ is of the form  $\left ( \begin{array}{cc} e^{t_k} & 0 \\
D_{k} & e^{2 t_k} \\ \end{array} \right ) $:

Considering if necessary a subsequence of $(b_k)$, and $(b_k^{-1})$ instead of $(b_k)$, we will suppose that $t_k$ has a limit in $\R_+^* \cup \{  + \infty \}$. So, looking at  a new subsequence of $(b_k)$, there exist two sequences of  $\RR^+$, $\beta_k$  and $\mu_k$, tending to  $\mu_{\infty}$ and  $\beta_{\infty}$ in  $\RR^+ \cup \{ + \infty \}$ and  $\RR_+^*\cup \{ + \infty \}$ respectively, as well as  a sequence  $B_k$ (resp. $\tau_k$) in  $End(\lien^+)$ (resp. in $S_{\lien^+}$) converging to $B_{\infty} \not = 0$ in  $ End(\lien^+)$ (resp. $\tau_{\infty} \in S_{\lien^+}$), such that $\rho(b_k)= \beta_kB_k + \mu_k\tau_k$, and such  that we are in one of the three following cases :

$\bullet $ {\it First case : $\lim_{k \to + \infty}\frac{\beta_k}{\mu_k}=0$}.

We set:
\[ F=\{ [x,u) \in {\bf [\dot Q]} \ | \ u  = - \tau_{\infty}  \}    \]

This is a closed subset of ${\bf [\dot Q]}$ (which is empty if $\tau_{\infty} \not \in S_{\liez^+}$). Let $[x_k,u_k)$ be a sequence of  ${\bf \Omega^+}={\bf [\dot Q]} \setminus F$ converging to $[x_{\infty},u_{\infty}) 
\in {\bf \Omega^+}$. One has $\rho(b_k)[x_k,u_k)=[x_k^{\prime},u_k)$, with $x_k^{\prime}= \mu_k(\frac{\beta_k}{\mu_k}B_k.x_k + \tau_k)$.  Under our hypothesis $\mu_k \to + \infty$ and thus  $\lim_{k \to + \infty} \frac{x_k^{\prime}}{||x_k^{\prime}||}= \tau_{\infty}$. Since $u_{\infty} \not = -\tau_{\infty}$, Lemma \ref{outil} applies and we get that $\lim_{k \to + \infty}\rho(b_k)[x_k,u_k)=[o]$.\\

If $x \in \liez^-$, we observe that $[\theta x,\frac{\theta x}{||x||})$ and  $[-\theta x, \frac{- \theta  x}{||x||})$ can't be in $F$ simultaneously. This proves that $s({\bf \Omega}^+)$ contains some elements of $\liez^-$ , which are arbitrarly close to $0_{\lieg}$. This imply the point $(ii)$ in this case.

$\bullet $ {\it Second case : $\lim_{k \to + \infty}\frac{\beta_k}{\mu_k}=+ \infty$}.

Let us set:

\[  F = \{  [x,u) \in {\bf [\dot Q]} \ | \ B_{\infty}.x \in \RR^-.u \}  \]

This is a closed subset of  ${\bf [\dot Q]}$. We first observe that in this case, we have necessarily $\beta_{\infty}=+ \infty$. If it were not the case, we should have $\mu_{\infty}=0$, and $\beta_{\infty} \in \RR_+^*$. In this case $L_k$ would converge to $L_{\infty} \in GL(\lieg)$, and the sequence  $(b_k)$ would be bounded. 

 Let  $[x_k,u_k)$ be a sequence of  ${\bf \Omega^+}={\bf [\dot Q]} \setminus F$ converging to $[x_{\infty},u_{\infty}) \in {\bf \Omega^+}$. One has  $\rho(b_k)[x_k,u_k)=[x_k^{\prime},u_k)$, with $x_k^{\prime} = \beta_k(B_k.x_k + \frac{\mu_k}{\beta_k} \tau_k)$, and  $\lim_{k \to + \infty}\frac{x_k^{\prime}}{|| x_k^{\prime} ||}= \frac{B_{\infty}.x_{\infty}}{|| B_{\infty}.x_{\infty} ||}$. By definition of  $F$, $B_{\infty}x_{\infty} \not \in \RR_-.u_{\infty}$. Lemma \ref{outil} applies, and we get  $\lim_{k \to + \infty} \rho(b_k)[x_k,u_k)=[o]$.\\

As we already noticed, $\beta_{\infty}=+ \infty$. In a basis compatible with the  graduation $\lien^+=\lien_1^+ \oplus \liez^+$, the matrix $B_{\infty}$ has either the form $\left ( \begin{array}{cc} 0 & 0 \\
A_{\infty} & 0 \\ \end{array} \right ) $, with $A_{\infty} \not = 0$ (this  possibility occures only if $\lieg \not = \mathfrak{so}(1,n), \mathfrak{su}(1,1)$ or $\mathfrak{sp}(1,1)$), or  of the form $\left ( \begin{array}{cc} 0 & 0 \\
C_{\infty} & Id_{\liez^+} \\ \end{array} \right ) $ (we can have $C_{\infty}=0$). 

In this last  case, if $x \in \liez^- \setminus \{ 0_{\lieg} \}$, then $\theta x \not \in {\RR}^-. \frac{\theta x}{||x||}$. As a consequence, $\liez^- \setminus \{ 0_{\lieg} \} \subset s({\bf \Omega}^+)$ and the point $(ii)$ is proved in this case. 

In the  first case, we pick a sequence  $y_n$ in $\lien_1^+$ which tends to   $0_{\lieg}$,  and such that  for every $n \in \NN$, $A_{\infty}.y_n \not \in \RR_-.\frac{\theta x}{|| x ||}$. Then for $x \in \liez^- \setminus \{ 0_{\lieg} \} $, $s([\theta x + y_n, \frac{\theta x}{|| x ||}))$ is a sequence of  $s({\bf \Omega}^+)$ which tends to  $x$, what proves the point $(ii)$.

$\bullet $ {\it Third case :  $\lim_{k \to + \infty}\frac{\beta_k}{\mu_k}= \alpha$, with  $\alpha \in \RR_+^*$}.

Let us observe, as in the previous case, that since  $(b_k)$ tends to infinity, we have necessarily  $\beta_{\infty}=\mu_{\infty}=+ \infty$ in this case.

We set:
\[  F = \{  [x,u) \in {\bf [\dot Q]} \ | \  \alpha B_{\infty}.x + \tau_{\infty} \in \RR^-.u \}  \]

 This is a closed subset of  ${\bf [\dot Q]}$. Let  $[x_k,u_k)$ be a sequence of ${\bf \Omega^+}={\bf [\dot Q]} \setminus F$ converging to $[x_
{\infty},u_{\infty}) \in {\bf \Omega^+} $. One has $\rho(b_k)[x_k,u_k)=[x_k^{\prime},u_k)$, and   $x_k^{\prime} = \mu_k(\frac{\beta_k}{\mu_k}B_k.x_k + \tau
_k)$. Since  $[x_{\infty},u_{\infty}) \in {\bf [\dot Q]} \setminus F$, we get $\alpha B_{\infty}.x_{\infty}+\tau_{\infty} \not =0$, so that $\lim_{k \to + \infty}
\frac{x_k^{\prime}}{|| x_k^{\prime} ||}=\frac{\alpha B_{\infty}.x_{\infty}+\tau_{\infty}}{||\alpha B_{\infty}.x_{\infty}+\tau_{\infty}||} $. Since $\alpha B_{\infty}.x_{\infty}+
\tau_{\infty} \not \in \RR_-.u_{\infty} $, we can once again apply Lemma \ref{outil}, and we conclude:  $\lim_{k \to + \infty} \rho(b_k)[x_k,u_k)=[o]$.\\

Since we are still in the case $\beta_{\infty}=+ \infty$, the possible matrices for $B_{\infty}$ are the same as in the previous case. One thus checks that when $x \in \liez^-$ is very close to  $0_{\lieg}$, $[\theta x, \frac{\theta x}{||x||})$ and $[- \theta x , - \frac{\theta x}{||x||})$ can't be in  $F$ simultaneously, what proves point $(ii)$.

\end{preuve}




\section{ Cartan  geometries}
\label{cartan}

Here, we just indroduce the basic material on Cartan's geometries. For more details, and the proofs of some of the lemmas and propositions below, we refer to the very good reference  \cite{sharpe}. In all this section, we consider a Cartan geometry $(M,B, \omega)$, modelled on $\X=G/P$ (for the following generalities, we don't make special asumptions, so that  $G$ is any Lie group, and $P$ a closed subgroup of  $G$).

\subsection{Parallelism and Riemannian metric on  $B$
}
\label{metrique}
Let $(E_1,....,E_s)$ be a basis of the Lie algebra $\lieg$, which is orthonormal for the metric  $< \ >_{\lieg}$. At every point $p $ of $B$, we define a frame ${\cal R}(p)=({E}_1^{\dagger}(p),...,{E}_s^{\dagger}(p))$ by  ${E}_k^{\dagger}(p)=\omega_p^{-1}(E_k)$. This frame field yields a parallelism $\cal R$ on $B$.  For any vector   $\xi  =  \Sigma \lambda_k E_k$ in $\lieg$, we define a vector field $\xi^{\dagger}$ on $B$ by $\xi^{\dagger}(p)=\Sigma \lambda_k E_k^{\dagger}(p)$. Such vector fields with constant coordinates with respect to ${\cal R}$ are called {\it parallel}. 

\begin{notation}
\label{parallele}
We will adopt the notation $\xi_p$ instead of  $\xi^{\dagger}(p)$. Also, if ${\bf \Lambda}$ is a subset of  $\lieg$, it determines a field of "parallel subsets" of  $TB$, defined by $\Lambda_p=\omega_p^{-1}({\bf \Lambda})$, for every   $p \in B$. For a sequence ${\bf \Lambda}_k$ of subsets, we will write ${\Lambda_{p,k}}=\omega_p^{-1}({\bf \Lambda}_k)$.

If $\lieh \subset \lieg$ is a subalgebra,  $H_{p,\lieh}$ denotes $\omega_p^{-1}(\lieh)$.
\end{notation}

\begin{definition}

We call  $ \rho$ the Riemannian metric on  $B$ for which ${\cal R}(p)$ is an orthonormal frame for all  $p \in B$.
\end{definition}

Let us notice that for all  $p \in B$,  $\rho_p=\omega_p^* (< \ >_{\lieg})$.


\subsection{Curvature}
The Maurer-Cartan form on the Lie group $G$ satisfies the so called  {\it  structure equation} : $d \omega_{G} + \frac{1}{2}[\omega_G,\omega_G]=0$. Nevertheless, for a given Cartan connection $\omega$, the 2-form  $\Omega=d\omega +  \frac{1}{2}[\omega, \omega]$ needs not to be trivial. One calls $\Omega$ {\it the curvature} of the  connection $\omega$. Here are some fundamental properties of the curvature form :

$(i)$ For every $b \in P$, $(R_b)^* \Omega= Ad(b^{-1}) \Omega$.

$(ii)$ $\Omega(X,Y)=0$  if $X$ or $Y$ is a vertical vector (i.e tangent to the fibers of $B \to M$).

A Cartan geometry $(M, B,{\omega})$ is said to be {\it flat } if $\Omega$ vanishes identically on  $B$. 

Let us notice that the equivariance property  $(i)$ implies that the vanishing of the curvature at a point of a fiber implies its vanishing on the whole fiber. In the following, we will sometimes say abusively that $\omega$ vanishes at $x \in M$, meaning that $\omega$ vanishes on the fiber over $x$.

\begin{exemple}
Let $U$ be an open subset of $\X$,  stable under the action of a discrete subgroup $\Gamma \subset G$. If the action of $\Gamma$ on $U$ is free and properly discontinuous, then the manifold  $M=\Gamma \backslash U$ is naturally endowed with a flat Cartan geometry, inherited  from that of  $\X$. The bundle  $B$ is, in this case, the quotient  $\Gamma \backslash G$ and the Cartan connection is $\overline \omega_G$, the $1$-form induced by  $\omega_G$ on $\Gamma \backslash G$.

\end{exemple}

\subsubsection{Regularity}

We now say a few words on the asumption made on the connection $\omega$ in Theorem \ref{theoreme1}. Let us suppose that $(M,B,\omega)$ is modelled on a space $\X=G/P$, such that the Lie algebra $\lieg$ is endowed with a $k$-grading $\lieg=\lieg_{-k} \oplus ... \oplus {\lieg_0} \oplus ... \oplus \lieg_k$, $k \in \NN^*$, and the Lie algebra of $P$ is $\liep = \lieg_0 \oplus ... \oplus \lieg_k$. Notice that for all $-k \leq i \leq j \leq k$, $[\lieg_i,\lieg_j] \subset \lieg_{i+j}$. We then have an $Ad(P)$-invariant filtration $\lieg^{-k}=\lieg \supset \lieg^{-k+1} \supset ... \supset \lieg^{k}=\lieg_{k}$, putting $\lieg^{i}=\lieg_i \oplus \lieg_{i+1} \oplus ... \oplus \lieg_k$. 

\begin{definition}[Regularity]
The Cartan connection $\omega$ is said to be regular, if for any $\xi \in \lieg^{i}$ and $\zeta \in \lieg^{j}$, $i,j<0$, and any $p \in B$,  $\Omega_p(\xi_p, \zeta_p) \in \lieg^{i+j+1}$.

\end{definition}

This notion of regularity is important, since in most of the cases where the equivalence problem is solved, the "canonical" Cartan connection is regular.

Let us see more closely what is going on in the case where $\lieg$ is a rank one simple Lie algebra.

$\bullet$ If the roots of $\lieg$ are $\{ -\alpha, \alpha \}$, then the regularity condition just means that the curvature $\Omega$ takes its values in $\liep$. One says in this case that $\omega$ is {\it torsionfree}.

$\bullet$ If the roots of  $\lieg$ are $\{  -2 \alpha, - \alpha, \alpha, 2 \alpha\}$, then the regularity condition  means that for any pair $\xi, \zeta$ in $\lieg^{-1}=\lien_1^- \oplus \liep$, then $\Omega_p(\xi_p,\zeta_p) \in \lieg^{-1}$ (for one, and hence for any $p \in B$).

\subsection{Developping curves }
\label{developpante}
One of the fundamental properties of a Cartan connection, is that it establishes a link between parametrized curves of $M$ passing through a given point, and curves of the model space $\X$.

\begin{definition}
Let $N$ be a manifold, $q$ a point of  $N$ and  $I \subset \RR$  an open interval containing $0$. We define $C^1(N,I,q)$ as the space of $C^1$ maps $\gamma : I \to N$, such that $\gamma(0)=q$.
\end{definition}

We now recall how to develop curves of $M$ into curves on $\X$ thanks to the Cartan connection. The proofs can be found in \cite{sharpe} (Lemma 4.12 p 208).

\begin{lemme}
\label{developpe}
Let $p$ be a point of  $B$, and  $I$ an open interval containing $0$. Then there exists a unique map ${\hat {\cal D}}_p : C^1(B,I,p) \to C^1(G,I,e))$ satisfying:

$(i)$ For every curve $\hat \gamma$ of $C^1(B,I,p)$, the curve  $\alpha = {\hat{\cal D}}_p(\hat \gamma)$ satisfies for every $t \in I$: $\omega({\hat \gamma}^{\prime}(t))=\omega_G(\alpha^{\prime}(t))$. 

$(ii)$ If $a \in C^1(P,I,a_0)$, and if we define  $R_a \hat \gamma$ to be the curve defined by  $(R_a \hat \gamma)(t)=R_{a(t)} \hat \gamma(t)$, then ${\cal D}_{p.a_0}(R_a \hat \gamma)=a_0^{-1}.(R_{a}\alpha$). 
\end{lemme}

\begin{corollaire}
\label{descente}
Let $\gamma \in C^1(M,I,x)$  and $\hat \gamma_1 : I \to B$, ${\hat \gamma}_2 : I \to B$ be two lifts of   $\gamma$ in $C^1(B,I,p)$. Then  ${\hat{\cal D}}_{p}(\hat \gamma_1)$ and ${\hat{\cal D}}_{p}({\hat \gamma}_2)$ project on the same curve of ${\bf X}$.
\end{corollaire}

\begin{definition} [Developping map] The corollary implies that for all $x \in M$, and $\hx \in B$ over $x$, there is a well defined map from  $C^1(M,I,x)$ to  $C^1({\bf X},I,o)$. This map is called {\it developping map} at $x$,  and denoted by  ${\cal D}_x^{\hx}$ . If $\gamma \in C^1(M,I,x)$, and if  $[\gamma]$ is the  associated geometric segment, we will write ${\cal D}_x^{\hx}([\gamma])$ instead of  $[{\cal D}_x^{\hx}(\gamma)]$.
\end{definition}



\subsubsection{Development for flat structures}
\label{plat}

We refer once again to  \cite{sharpe} for the details concerning  this section.

The vanishing of the curvature is the only obstruction for the bundle $(B, \omega)$ to be locally isomorphic to $(G, \omega_G)$.

Let $\tilde M$ be the universal cover of $M$, and $r : \tilde M \to M$ the covering map. There exists a covering $\tilde r : \tilde B \to B$, which is a  principal $P$-bundle $\tilde \pi : \tilde B \to \tilde M$, such that the following diagram is commutative:

$$\begin{array}{ccc} 
\tilde B & \stackrel{\tilde r}{\rightarrow} & B \\
 \downarrow & &\downarrow  \\
\tilde M & \stackrel{r}{\rightarrow} & M \\

\end{array}$$

The form  ${\tilde \omega}={\tilde r}^* \omega$ is a Cartan connection on $\tilde B$.

\begin{theoreme}
If the curvature of $\omega$ vanishes identically, then there exists a bundle morphism  $\tilde \delta :  (\tilde B, \omega) \to (G, \omega_G)$, which is an immersion satisfying ${\tilde \delta}^* \omega_G={\tilde \omega}$. The map $\tilde \delta$ induces an  immersion $\delta : \tilde M \to \X$. 
\end{theoreme}

In that case, the manifold  $M$ is said to be endowed with a  $(G, \X)-$structure, and the map $\delta$ is known as a {\it developping map} of the structure.
Let $(M, B, \omega)$ be a flat Cartan geometry, $x$ a  point of $M$ and $\tilde x$ a  point of $\tilde M$ over $x$. We can compose  $\tilde \delta$ by an element of $G$ in order to get $\delta(\tilde x)=o$. Let $\gamma \in C^1(M,I,x)$.  Then there is a  unique $\tilde \gamma \in C^1(\tilde M,I,\tilde x)$,  such that $r(\tilde \gamma(t))=\gamma(t)$ for every $t \in I$. If $\hx$ is the point over $x$ such that $\tilde \delta (\hx)=e$, then ${\hat{\cD}_x^{\hx}}(\tilde \gamma)=\delta \circ \tilde \gamma$.

\section{The geodesics of a rank one parabolic geometry}
\label{geodesique}

\subsection{Exponential maps}

For  $p \in B$ and   $\xi \in T_pB$, we define  $U=\omega_p(\xi)$, and  $\phi_\xi^t$ the local flow associated to the parallel vector field $U^{\dagger}$. There is a neighbourhood $W_p^B$ of $0_p$ in $T_pB$, starshaped with respect to $0_p$, such that for every $\xi \in W_p^B$, $t \mapsto \phi_\xi^t(p)$ is defined on $[0,1]$. 

\begin{notation}
\label{nono}
We will denote by  $W^{\cB}$ the open subset of  $T{\cB}$  defined by $\bigcup_{p \in \cB}W_p^B$.\\
If $p \in B$,  ${\bf W}_{p}$ is the subset of $\lieg$ defined by  $\omega_p(W_p^B)$.
\end{notation}

\begin{definition}

We define an  {\it exponential map} $Exp : W^{\cB} \to \cB$, by $Exp_p(\xi)=\phi_\xi^1(p)$, for every $p \in \cB$, $\xi \in W_p^B$.

-If  $x = \pi(p)$ and  $\xi \in W_p^B$, we define  $Exp_x(\xi)=\pi \circ Exp_p(\xi)$.

- For $\xi \in W_p^B$, we call ${\hat \xi}^*$ (resp. ${\xi}^*$) the path of  $B$ (resp. of  $M$), parametrized by $[0,1]$, and  defined by  ${\hat \xi}^*(t)=Exp_{p}(t\xi)$ (resp. ${\xi}^*(t)=Exp_x(t \xi)$). The geometric support of ${\hat \xi}^*$ (resp. ${\xi}^*$) will be denoted by  $[{\hat \xi}]$ (resp. $[{\xi}]$).

- We will  write  $[x]$ instead of $[0_p]$.
\end{definition}

- If  $\Lambda_p \subset W_p^B$, we set $[\Lambda_p]=\bigcup_{\xi \in \Lambda_p}[\xi]$. Let us notice that $Exp_{\pi(p)}(\Lambda_p) \subset [\Lambda_p]$.

\subsection{Geodesics}

We now suppose that  $(M,B, \omega)$ is a Cartan geometry modelled on  $\X=\partial {\bf H}_{\KK}^d$.

\begin{definition} [Geodesics on $M$]
Let  $I$ be an interval containing $0$, $x$ a point  of $M$, $\hat x$ over $x$, and   $\gamma \in C^1(M,I,x)$. One says that  $\gamma$ is a parametrized geodesic of  $M$ if and only if ${\cD}_x^{\hat x}(\gamma)$ is a parametrized geodesic of  ${\bf X}$, as defined in section  \ref{geox}. This definition does not depend of the choice of  $\hx$ over $x$, by $Ad(P)$-invariance of ${\bf Q}$
\end{definition}

\begin{definition}
- For every $p \in B$, we define: 
\[ Q_p= \{\xi \in T_pB \ | \ \omega_p(\xi) \in {\bf Q}    \}   \]
and  ${\dot Q}_p=Q_p \setminus \{  0_p \}$

- $TQ=B \times_P {\bf Q}$ is the subbundle of $TB$, the fibers of which are the   $Q_p$'s.

- For every $p \in B$,  we call $Q_p^{reg}$ the {\it regular part } of $Q_p$, namely the set of $u \in {\dot Q}_p \cap W_p^B$such that the restriction of  $Exp_{\pi(p)}$ to  $Q_p$ is a submersion at $u$.
\end{definition}


\begin{lemme}
Let  $x$ be a point of  $ M$, and  $\hat x \in B$ over$x$. The geodesic segments parametrized by $[0,1]$, starting from $x$, are exatly the $\xi^*$'s, for  $\xi \in W_{\hat x}^B \cap Q_{\hx}$. 

\end{lemme}

\begin{preuve}

By the very definition of the parallel fields, we get the relation :

\begin{equation}
\label{commute}
{\cD}_x^{\hx}(\xi^*)=(\omega_{\hat x}(\xi))^*
\end{equation}

The lemma then follows immediately.
\end{preuve}

Now comes a very important lemma, which links the behaviour of a sequence of geodesic segments in $M$ (or sets of such sequences) with the behaviour of their developments:

\begin{lemme}
\label{collapse}

Let $\hx_k$ be a sequence of  $B$ converging to $\hx_{\infty} \in B$,  and  ${\bf \Lambda}_k$ a sequence of subsets of ${\bf Q} \cap {\bf W}_{\hx_k}$. Then if  $\lim_{k \to + \infty}[{\bf \Lambda}_k]=[o]$, we also have $\lim_{k \to + \infty} [{\Lambda_{\hx_k,k}}]=[x_{\infty}]$ (where, as usual,   $\Lambda_{\hx_k,k}$ stands for  $\omega_{\hx_k}^{-1}({\bf \Lambda}_k)$). 
\end{lemme}

\begin{preuve}
For every $p \in B$, the space $H_{p, \lien^-}$ (see the notation  \ref{parallele}) is called {\it horizontal space } at $p$ (notice that this distrubution is not invariant by the action of  $P$ on $B$, unlike  the case of Ehresmann's connections).

Let $x \in M$, $\hx \in B$ over $x$, and $\gamma \in C^1(M,I,x)$,  such that   $\cD_x^{\hx}(\gamma) \subset \Omega_{o}$. Then, there exists a unique horizontal lift   $\gamma$ in $C^1(B,I,\hat x)$ ( i.e with horizontal tangent vector for all $t \in I$). Indeed, let $\hat \gamma$ be any lift of  $\gamma$ in $C^1(B,I, \hat x)$, and $\hat \alpha = {\hat{\cD}}_{\hx}(\hat \gamma)$. Since  $N^-$ is transverse to any fiber of $G \to \X$ over $\Omega_{o}$, there is  $a \in C^1(P,I,e)$ such that  $R_a\hat \alpha$ is a curve in $C^1(N^-,I,e)$. Thanks to Lemma  \ref{developpe}, point $(ii)$, the curve  $R_a \hat \gamma$ is also horizontal.

Now, let us go back to the proof of lemma  \ref{collapse}. Suppose first that ${\bf \Lambda}_k$ is just a sequence of points ${ \xi}_k$ in ${\bf Q} \cap {\bf W}_{\hx_k}$. By the previous remark, there is an horizontal curve  $\hat \gamma_k \in C^1(B,I, \hat x_k )$ over  $\xi_{\hx_k,k}^*$. The curve $\hat \alpha_k = {\hat{\cD}}_{\hx_k}(\hat \gamma_k)$ is in  $C^1(N^-,I,e)$ and projects on ${\bf \xi}_k^*$.  Now, we observe that if $(a_{1,k}(t),...,a_{n,k}(t))$ are the coordinates of ${\hat \gamma}_k^{\prime}(t)$ with respect to the parallelism ${\cal R}$, then the coordinates of  ${\hat \alpha}_k^{\prime}(t)$ with respect to the parallelism defined by $E_1,...,E_n$ on $TG$ are also $(a_{1,k}(t),...,a_{n,k}(t))$. We thus get that $L_{\rho}({\hat \gamma}_k)=L_{\rho_G}({\hat \alpha}_k)=L^-(\xi_k^*)$. Thanks to the hypothesis of the lemma, and property $(P_3)$ of lemma  \ref{pro}, we get  : $\lim_{k \to + \infty} L_{\rho}(\hat \gamma_k)= 0$.  So, for $k$ sufficiently large, $[{\hat \gamma}_k]$ is included in any  $\rho$-ball of arbitrary small radius, and centered at  $\hx_{\infty}$. Projecting, $[{\xi_{\hx_k,k}}]$ is included, for  $k$ large enough, in any neighbourhood of $x_{\infty}$, what concludes the proof. 

Now, in the general case, $[\Lambda_{\hx_k,k}]$ tends to $[x_{\infty}]$ if and only if for any sequence $\sigma_k$ of $\Lambda_{\hx_k,k}$, $[\sigma_k]$ tends to $[x_{\infty}]$, so that the proof of the general case reduces to that of the previous situation.

\end{preuve}

We end the section with a technical lemma:

\begin{lemme}
\label{submersion}
Let  $x$ be a point of $  M$ and  $\hat x$ a point of  $B$ over $x$. There is an open subset $U \subset H_{\hat x, \liez^-}$,  containing  $0_{\hat x}$, such that   $U \setminus \{  0_{\hx} \} \subset Q_{\hx}^{reg}$.
\end{lemme}

\begin{preuve}
Let us first remark that for all $\xi \in \liez^-$, $\lien^-  \subset T_\xi{\bf {\dot Q}}$. This is quite obvious when   $\liez^-=\lien^-$. In the other cases, let us fix  $\zeta_1,....,\zeta_s$ a basis of  $\lien_1^+$.  The flow $exp_G(t\xi)$  commutes  with none of the $Exp_G(t\zeta)$'s, for any $\zeta$ in $\lien^-$. Indeed, the first flow acts on  $\X$ with a unique fixed point which is $\nu$. The second ones act on  $\X$ also with a unique fixed point,  which is  $o$.  As a consequence, $[\xi,\zeta_1],...,[\xi,\zeta_s]$ is a free family of vectors, which turns out to be a basis of  $\lien_1^-$ (indeed, we are in the case where  $\liez^-=\lieg_{-2 \alpha}$ and  $\lien_1^+=\lieg_{+ \alpha}$, see section \ref{algebriques}. So,  $ [\liez^-,\lien_1^+] \subset \lien_1^-$). On the other hand, since for every  $1 \leq i \leq n$, $[\xi,\zeta_i]=\frac{d}{dt}_{|t=0}Ad(Exp_G(t \zeta_i)).\xi$, we get $[\xi,\zeta_i]   \in T_{\xi}{\bf \dot Q}$. Finally, $\lien_1^-   \subset T_{\xi}{\bf \dot Q}$.  By definition of ${\bf Q}$, $\liez^- \subset T_{\xi}{\bf \dot Q}$. We then get $\lien^- \subset T_\xi{\bf {\dot Q}}$.

 Let us choose an open subset $V_{\hat x}$ containing $0_{\hat x}$,  such that   $Exp_{\hat x}$ is a diffeomorphism from $V_{\hat x}$ onto its image. The pullback in $V_{\hat x}$ by $Exp_{\hat x}$, of the fibers of $B \to M$,  yields a foliation ${\mathcal F}$ of $V_{\hat x}$. The leaf  ${\cal F}_{0_{\hat x}}$ passing through  $0_{\hat x}$ is nothing else than  $\omega_{\hat x}^{-1}(\liep) \cap V_{\hat x}$. In particular, it is transverse to  $H_{\hat x, \liez^-}$. Thus, there is a neighbourhood $U$ of  $0_{\hat x}$  in $H_{\hat x, \liez^-}$,  such that for every $u \in U$, ${\mathcal F}_u$ is  transverse to  $H_{\hat x, \lien^-}$. Now, by what we said before, $H_{\hx,\lien^-} \subset T_u{{\dot Q}_{\hx}}$ for every  $u \in U$. The differential of  $Exp_{x}$ at $u$, when restricted to  $T_u{{\dot Q}_{\hx}}$ is then a  surjection onto $T_{Exp_x(u)}M$.

\end{preuve}

\section{Automorphisms of a  Cartan Parabolic geometry}

\subsection{Preliminary remarks}
In this subsection, we consider here a general Cartan geometry. In the introduction of the article, we defined the group $Aut(B, \omega)$, as the group of  diffeomorphisms $h$ in $B$ such that  $h^* \omega= \omega$. Any element of  $Aut(B, \omega)$ preserves the parallelism ${\cal R}$, so that  $Aut(B, \omega)$ is a closed subgroup of  $Iso(B, \rho)$. This last group is closed in $Homeo(B)$, by classical properties of the group of isometries of a Riemannian manifold (se  for example \cite{kobayashi}, \cite{sternberg}).

\subsection{Holonomy sequences}
\label{holonomie}

\begin{definition}[Holonomy datas]
Let  $(f_k,x_k,x_{\infty},y_{\infty})$ be a stable data of $Aut(M,\omega)$. We say that 
$(b_k, \hx_k,\hx_{\infty},\hy_{\infty})$ is an associated holonomy data, if there is a lift ${\hat f}_k$ of $f_k$ in $Aut(B, \omega)$, 
 a sequence $\hx_k$  of $B$ over  $x_k$, tending to  $\hx_{\infty}$, and a sequence
 $\hy_{\infty} \in B$  over  $y_{\infty}$, so that  $\hy_k=R_{b_k}^{-1} \circ \hat f_k(\hx_k)$ tends to $\hy_{\infty}$.

\end{definition}

To any stable data $(f_k,x_k,x_{\infty},y_{\infty})$ of $Aut(M, {\cal S})$, it is possible to associate an holonomy data $(b_k, \hx_k,\hx_{\infty},\hy_{\infty})$.  To see this, let us fix two open subsets $U$ and $V$ around  $x_{\infty}$ and  $y_{\infty}$ respectively. We suppose $U$ and $V$ so small that there are two continuous sections $s_U:U \to B$ and $s_V:V \to B$. Let ${\hat f}_k$ be a lift of $f_k$ in $Aut(B, \omega)$. We then define   $s_U(x_k)=\hat x_k$, $k \in \NN \cup \{ \infty \} $ (resp. $s_V(f_k(x_k))=\hat y_k$, $k \in \NN \cup \{ \infty \} $). Since  $\hat f_k(\hx_k)$ is a sequence over $f_k(x_k)$, there is a unique sequence  $(b_k)$ of $P$, such that for all $k \in \NN$, we have $R_{b_k^{-1}} \circ \hat f_k(\hx_k)=\hy_k$.

\begin{proposition}[equivariance]
\label{equivariance}
Let $(f_k,x_k,x_{\infty},y_{\infty})$ be a stable data of  $Aut(M,\omega)$, and $(b_k,\hx_k,\hx_{\infty},\hy_{\infty})$ an associated holonomy.  For every sequence $\xi_k \in {\bf W}_{\hx_k}$, we  define $\zeta_k=Ad(b_k).\xi_k$. We then have the equivariance property:

\begin{equation}
\label{commute3}
f_k(\xi_{\hx_k,k}^*)=\zeta_{\hy_k,k}^* 
\end{equation} 
\end{proposition}

\begin{preuve}
Let ${\hat f}_k$ be a lift of $f_k$ defining the holonomy sequence $b_k$, and put $\phi_k=R_{b_k^{-1}} \circ {\hat f}_k$. By definition $\phi_k(\hx_k)=\hy_k$ and $\lim_{k \to +\infty}\hy_k=\hy_{\infty}$. We observe that $(\hat \phi_k)_{*}{\xi^{\dagger}}=(Ad(b_k)\xi)^{\dagger}$. We get  firstly that if  $\xi_k \in {\bf W}_{\hx_k}$, then  $\zeta_k \in {\bf W}_{\hy_k}$, and secondly that $\hat \phi_k({\hat \xi}_{\hx_k,k}^*)={\hat \zeta}_{\hy_k,k}^*$. Projecting on $M$, this leads to the relation  (\ref{commute3}).

\end{preuve}

\subsection{The rank one case}
Now, we consider a Cartan geometry $(M, B, \omega)$,  modelled on some space  ${\X}= \partial {\bf H}_{\KK}^d$.

\begin{theoreme} 
\label{automorphisme}
\ \\
$(i)$ The group $Aut(M, \omega)$ is closed in  $Homeo(M)$.  \\

Let  $(f_k,x_k,x_{\infty},y_{\infty})$  be a stable data of $Aut(M,\omega)$, and  $(b_k,\hx_k,\hx_{\infty},\hy_{\infty})$ an associated holonomy. Then:

$(ii)$ $(f_k)$ is bounded in  $Aut(M, \omega)$ if and only if  $(b_k)$ is bounded in  $G$.

$(iii)$ If the action of  $(f_k)$ is equicontinuous at  $x_{\infty}$, then the action of $(b_k)$ is equicontinuous at $o$.
\end{theoreme}

\begin{preuve}





Let us begin with the proof of the  point $(iii)$.

Suppose on the contrary that  $(b_k)$ does not act equicontinuously at  $o$. Looking if necessary at a subsequence of $(b_k)$, one can  choose a sequence  $(z_k)$ in $\X$, converging to $o$, and such that  $w_k=b_k.z_k$ tends to $w_{\infty} \not = o$. By the property  $(P_1)$ of lemma \ref{pro}, there exists, for all  $k$, a geodesic segment  $\alpha_k=[z_k,u_k)$ linking $o$ to $z_k$ (we still identify ${\bf [Q]}$ with $\lien^+ \times {S}_{\liez^+}$). We choose the  $\alpha_k$'s such that $\lim_{k \to + \infty} \alpha_k = [o]$. Looking at a subsequence if necessary, we suppose that $\beta_k=b_k.\alpha_k$ tends to $\beta_{\infty}=[w_{\infty},u_{\infty})$. We put  $\xi_k=s(\alpha_k)$ and  $\zeta_k=s(\beta_k)$. We also set  $\zeta_{\infty} =s(\beta_{\infty})$. By continuity of $s$,  $\lim_{k \to + \infty}\zeta_k=\zeta_{\infty}$.   We can find a real $\lambda \in ]0,1]$ such that  $\lambda \zeta_k \in {\bf W}_{\hy_k}$ for all $k \in \NN \cup \{  \infty \}$. So, replacing if necessary $\zeta_k$ by $\lambda \zeta_k$, and $\xi_k$ by $Ad(b_k^{-1})(\lambda \zeta_k)$, we will suppose that  $\xi_k \in  {\bf W}_{\hx_k}$ for all $k \in \NN$, and $\zeta_k \in {\bf W}_{\hy_k}$ for all $k \in \NN \cup \{  \infty \}$.   

From Proposition \ref{equivariance}, we infer that for all $k \in \NN $, and $t \in [0,1]$ : $f_k(\xi_{\hx_k,k}^*(t))=\zeta_{\hy_k,k}^*(t)$. Since $\zeta_{\hy_{\infty},\infty} \not = 0$, there exists $t_0 \in [0,1]$ such that $\zeta_{\hy_{\infty},\infty}^*(t_0) \not = y_{\infty}$.  So, $\lim_{k \to +\infty}f_k(\xi_{\hx_k,k}^*(t_0)) \not = y_{\infty}$. Nevertheless,  $\xi_{\hx_k,k}^*(t_0)$ tends to $x_{\infty}$. Indeed, $\lim_{k \to + \infty}[\xi_k]=[o]$ and Lemma \ref{collapse} implies that $\lim_{k \to + \infty}[\xi_{\hx_k,k}]=[x_{\infty}]$. We conclude that the action of  $(f_k)$ is not equicontinuous at $x_{\infty}$.\\

We can now prove the point $(i)$. We suppose that $(f_k)$ is a sequence of $Aut(M, \omega)$ converging to $f_{\infty} \in Homeo(M)$. Then, for all point $x \in M$, if $( x_k)$ is the constant sequence equal to $x$, and $y_{\infty}=f_{\infty}(x)$, the quadruple $(f_k, x_k,x,y_{\infty})$ is a stable data.  Let us denote by $(b_k, \hx_k,\hx,\hy_{\infty})$ an associated holonomy data. We prove :

\begin{lemme}
The sequence $(b_k)$ is bounded in $G$.

\end{lemme}

\begin{preuve}
If $(b_k)$ were unbouded, there would be a subsequence of $(b_k)$ admitting a North-South dynamic with poles $o^+$ et $o^-$. Since $o$ is fixed by $(b_k)$, we have necessarily $o=o^+$ ou $o=o^-$. In the second case, the action of $(b_k)$is not equicontinuous at $o$ (see Corollary  \ref{kak}), and as a consequence of point $(iii)$,  the action of $(f_k)$ should not be equicontinuous at  $x$. This is in contradiction with the fact that $(f_k)$ has a limit in $Homeo(M)$. Thus, we are in the case  $o=o^+$. Let $\alpha=[\xi]$ be a segment of  ${\bf [Q]} \setminus [o]$, such that $\xi \in {\bf W}_{\hx}$, and $\alpha \subset \Omega_o$. Let us remark that $\zeta_k=Ad(b_k).\xi \in {\bf W}_{\hy_k}$ for all $k \in \NN$.  Since $o=o^+$, we have $\lim_{k \to + \infty} b_k.\alpha = [o]$. By Lemma \ref{collapse}, we obtain  $\lim_{k \to + \infty}[\zeta_{\hy_k,k}]=[y_{\infty}]$. But Proposition  \ref{equivariance} gives that  $f_k([\xi_{\hx_k}])=[\zeta_{\hy_k,k}]$ for all $k \in \NN$, which leads to  $f_{\infty}([\xi_{\hx}	])=[y_{\infty}]$. Since $[\xi_{\hx}] \not = [x]$, we get a contradiction with the fact that  $f_{\infty} \in Homeo(M)$.

\end{preuve}

We take again a subsequence of $(f_k)$, so that $(b_k)$ tends to $b_{\infty} \in P$. 
 This means that $\hat f_k(\hat x_k)$ tends to  $R_{b_{\infty}}.\hat y_{\infty}$. But $(\hat f_k)$ is a sequence of isometries for the riemannian metric $\rho$. Since the action of the isometry group of a riemannian manifold is proper, $\hat f_k$ is relatively compact in  $Iso(B,\rho)$. The group $Aut(B,\omega)$ being closed in  $Iso(B,\rho)$, we can suppose,  looking once again at a subsequence, that there exist  $ \hat f_{\infty} \in Aut(B, \omega)$ such that $\hat f_k \to \hat f_{\infty}$. Finally, $\hat f_{\infty}$ is a lift of $f_{\infty}$ in $Aut(B, \omega)$, so that  $f_{\infty} \in Aut(M,\omega)$.\\
 
 The proof of point $(ii)$ goes in a way similar to that of point  $(i)$. 

\end{preuve}

\begin{remarque}
It is likely that this use of geodesics and holonomy can be generalized, to prove that the automorphisms group of any parabolic Cartan geometry on $M$ is closed in $Homeo(M)$. 
\end{remarque}
\section{Proof of Theorem \ref{theoreme1}}
\label{demonstrations}

We first make the proof of Theorem  \ref{theoreme1} in an easy case: the case of Kleinian manifolds (i.e manifolds $M$ which are a quotient of some open subset $\Omega \subset \X$ by a discrete subgroup  $\Gamma \subset G$). These manifolds are naturally endowed with a Cartan geometry modelled on $\X$, with $B= \Gamma \backslash \pi_X^{-1}(\Omega)$, and $\omega= \overline \omega_G$, the Cartan connection induced by  $\omega_G$ on this quotient.

\begin{lemme}
\label{bebe-obata}
Let $\Omega$ be a connected open subset of  $\X$, and $M= \Gamma \backslash \Omega$, where $\Gamma \subset G$ is a discrete subgroup acting freely properly discontinuously on  $\Omega$. 

If $Aut(M, \overline \omega_G)$ does not act properly on  $M$, then $\Gamma = \{  e \}$ and:

$\bullet$ either $\Omega= \X \setminus \{ \kappa \}$, for some $\kappa \in \X$.

$\bullet$ or $\Omega=\X$.

\end{lemme}

\begin{preuve}
The group $Aut(M, \overline \omega_G)$ is induced by elements of $G$ leaving $\Omega$ stable and normalizing $\Gamma$. If $Aut(M, \overline \omega_G)$ does not act properly on $M$, we can find a sequence $(h_k)$ of $G$,  normalizing  $\Gamma$, and acting nonproperly on $\Omega$. We can suppose that this sequence tends to infinity in $G$, and looking if necessary at $(h_k^{-1})$ instead of $(h_k)$, we can also asume that $o^- \in \Omega$. Let us call  $\pi_M$ the covering map  from $\Omega$ onto $M$, and let $U \subset \Omega$ a small open subset containing  $o^-$,  such that  $\pi_M$ maps $U$ diffeomorphically on its image. By Lemma \ref{nord-sud}, there is a sequence $(k_m)_{m \in \NN}$ such that  $h_{k_m}(U)$ is an increaszing sequence of open subsets, the union of which is the whole space $\X$ if $o^+=o^-$, and  $\X \setminus \{  o^+ \}$ if $o^+ \not = o^-$. We infer that $\Omega = \X$ or $\Omega = \X \setminus \{ o^+ \}$. Moreover, since  $h_k$ normalizes $\Gamma$ for all $k \in \NN$, $\pi_M$ has to be injective on each $h_{k_m}(U)$. Finally, $\pi_M$ is injective on a dense open set of $\Omega$, which implies $\Gamma = \{  e \}$.

\end{preuve}

We now work under the general hypothesis of Theorem  \ref{theoreme1}: $(M,B, \omega)$ is a Cartan geometry modelled on $\X= \partial {\bf H}_{\KK}^d$, $\omega$ is a regular connection, and  $Aut(M,\omega)$ does not act properly on $M$. 

\begin{proposition}
\label{P}
If $(f_k)$ is a sequence of  $Aut(M, \omega)$ which does not act properly on $M$, then there is an open subset  $U \subset M$ and a point $y_{\infty} \subset M$, such that  $\lim_{k \to + \infty }f_k(U)=y_{\infty}$.

\end{proposition}

\begin{preuve}
Since  $Aut(M,\omega)$ does not act properly on $M$, there is a stable data  $(f_k,x_k,x_{\infty}, y_{\infty})$, with $(f_k)$ a sequence of  $Aut(M, \omega)$ tending to infinity. Let $(b_k,\hat x_k, \hx_{\infty}, \hy_{\infty})$ be an associated holonomy.  By theorem  \ref{automorphisme}, the sequence $(b_k)$ tends to infinity in $P$. So, changing if necessary  $(f_k,x_k,x_{\infty}, y_{\infty})$ into  the stable data $(f_k^{-1}, y_k,y_{\infty}, x_{\infty})$, Proposition \ref{dynamique} ensures the existence of an open subset ${\bf \Omega}^+ \subset {\bf [\dot Q]}$, such that for any  compact subset $K \subset { \bf \Omega}^+$, $\lim_{k \to + \infty} b_k.K = [o]$. Let ${\bf U}^+=s({\bf \Omega}^+)$. 

As a consequence of Lemma \ref{submersion} and  point $(ii)$ of Proposition \ref{dynamique}, we can find $\xi \in {\bf U}^+$ such that $\xi_{\hx_{\infty}} \in \dot Q_{\hx_{\infty}}^{reg}$.   In fact, there exists $\epsilon >0$ small,  such that $B_{\xi_{\hx_{\infty}}}(\epsilon) \subset \dot Q_{\hx_{\infty}}^{reg}$, where $B_{\xi_{\hx_{\infty}}}(\epsilon)=\omega_{\hx_{\infty}}^{-1}(B_{\rho_G}(\xi,\epsilon) \cap {\bf U}^+)$.   This property still holds if one replaces $x_{\infty}$ by $x_k$,  for $k \geq k_0$ sufficiently large. In other words, if we put  $B_{\hx_k}(\epsilon)=\omega_{\hx_k}^{-1}(B_{\rho_G}(\xi, \epsilon) \cap {\bf U}^+)$, then, for $k \geq k_0$ ,  $U_k=Exp_{x_k}(B_{\hx_k}(\epsilon))$ is a sequence of open subsets of  $M$. Moreover, this sequence converges to $U_{\infty}=Exp_{x_{\infty}}(B_{\hx_{\infty}}(\epsilon))$. This ensures the existence of an open subset $U \subset M$ and of $k_1 \geq k_0$, such that  $U \subset \bigcap_{k \geq k_1}U_k$. 

Now, we suppose $\epsilon$ small enough, such that  $B_{\rho_G}(\xi,\epsilon) \cap {\bf U}^+$ has compact closure in ${\bf U}^+$.  Let us set ${\bf \Lambda}_k=Ad(b_k).(B_{\rho_G}(\xi,\epsilon) \cap {\bf U}^+)$.   We then have  $\lim_{k \to + \infty} [{\bf \Lambda}_k]=\lim_{k \to + \infty} b_k.[B_{\rho_G}(\xi,\epsilon) \cap {\bf U}^+]=[o]$. We also get, as a consequence of relation (\ref{commute3}) in Proposition \ref{equivariance}, that $[\Lambda_{\hx_k,k}]=f_k([B_{\hx_k}(\epsilon)])$, .   Lemma \ref{collapse} then gives  $\lim_{k \to + \infty} f_k([B_{\hx_k}(\epsilon)])= [y_{\infty}]$. Since $Exp_{x_k}(B_{\hx_k}(\epsilon)) \subset [B_{\hx_k}(\epsilon)]$, we get $U \subset [B_{\hx_k}(\epsilon)]$ for all $k \geq k_1$, and the proposition follows.


\end{preuve}

\begin{proposition}
\label{courbure}
Let $(f_k)$be a sequence of  $Aut(M,\omega)$ tending to infinity, and  $(f_k,x_k,x_{\infty},y_{\infty})$ a stable data on $M$. If the action of  $(f_k)$ is equicontinuous at  $x_{\infty}$, then the curvature of  $\omega$ vanishes at $x_{\infty}$.
\end{proposition}

\begin{preuve}
Let $(b_k,\hx_k,\hx_{\infty},\hy_{\infty})$ be an holonomy associated to   $(f_k,x_k,x_{\infty},y_{\infty})$. By Theorem  \ref{automorphisme}, the sequence  $(b_k)$ acts equicontinuously at  $o$.
By Corollary \ref{kak}, the attracting and repelling poles $o^+$ and  $o^-$ of $(b_k)$ satisfy $o^+=o$  and $o^- \not = o$, and we can write $b_k=l_{1,k}a^{t_k}l_{2,k}$, where $l_{1,k}$ and $l_{2,k}$ are two sequences of$P$, relatively compact in   $P$. We will suppose that  for $i \in \{ 1,2 \}$, $l_{i,k}$ converges to  $l_{i,\infty} \in P$. Since $o^+=o$, we have $\lim_{k \to + \infty}t_k=+ \infty$.

Let us fix $\eta_1,...,\eta_n$, a basis of  $\lien^-$, and we suppose that $\{ 1,...,n  \}=I_{-1} \cup I_{-2}$, with $\eta_i \in \lien_1^-$ for$i \in I_{-1}$ and $\eta_i \in \liez^-$ for $i \in I_{-2}$ (of course $I_{-1} = \emptyset$ when $\lien_1^-=\{ 0 \}$). For $i=1,...,n$ we set  $\xi_{i,k}=Ad(l_{2,k}^{-1}).\eta_i$ (resp. $\zeta_{i,k}=Ad(l_{1,k}).\eta_i$). For every $i=1,...,n$, we have  $\lim_{k \to + \infty} \xi_{i,k}=\xi_{i,\infty}$ (resp. $\lim_{k \to + \infty} \zeta_{i,k}=\zeta_{i,\infty}$), where  $\xi_{i,\infty}=Ad(l_{2,\infty}^{-1}).\eta_i$ (resp.  $\zeta_{i,\infty}=Ad(l_{1,\infty}^{-1}).\eta_i$).  In the following, we will denote $\lieh_k$ (resp. $\lieh_{\infty}$) the vector subspace of  $\lieg$ spanned by $\xi_{1,k},...,\xi_{n,k}$ (resp. $\xi_{1,\infty},...,\xi_{n,\infty}$). Since  $\lieh_k$ is nothing else than the image of  $\lien^-$ by $Ad(l_{2,k}^{-1})$, the $\lieh_k$'s (resp. $\lieh_{\infty}$) are Lie subalgebras of $\lieg$.

Let $({\hat f}_k)$ be a lift of $(f_k)$ associated to the holonomy $(b_k)$, i.e such that $R_{b_k^{-1}} \circ \hat f_k(\hx_k)=\hy_k$. Let us set  $\phi_k = R_{b_k^{-1}} \circ \hat f_k$. Then,  $(\phi_k)^*\Omega=Ad(b_k).\Omega$. We now distinguish two cases.\\

$a)$ The root system of $\lieg$ is $\Delta = \{ -\alpha, \alpha \}$. Then, for every $k \in \NN$, and $1 \leq i \leq n$, we have $Ad(b_k).\xi_{i,k}= e^{-t_k}\zeta_{i,k}$. We infer the  relation $D_{\hx_k}\phi_k(\xi_{\hx_k,i,k})=e^{-t_k}\zeta_{\hy_k,i,k}$, which implies for every $1 \leq i \leq j \leq n$ :
\[ Ad(b_k).\Omega_{\hat x_k}(\xi_{\hx_k,i,k},\xi_{\hx_k,j,k})= e^{-2t_k} \Omega_{\hat y_k}(\zeta_{\hy_k,i,k},\zeta_{\hy_k,j,k})  \]

Now, looking at the root space decomposition of $\lieg$, there is a $C>0$ such that for any $u \in \lieg$, $|| Ad(b_k).u|| \geq e^{-t_k}. || u ||$. We infer that $\Omega_{\hat x_{\infty}}(\xi_{\hx_{\infty},i,\infty},\xi_{\hx_{\infty},j,\infty}))= \lim_{k \to + \infty} \Omega_{\hat x_k}(\xi_{\hx_k,i,k},\xi_{\hx_k,j,k})=0$. The subspace ${\lieh}_{\infty}$ is a supplementary of $\liep$ in $\lieg$, so that $\Omega_{\hx_{\infty}}=0$.\\

$b)$ The root system of $\lieg$ is $\Delta =\{ -2 \alpha, - \alpha,  \alpha, +2 \alpha \} $.
Then, $Ad(b_k).\xi_{i,k}= e^{-t_k}\zeta_{i,k}$ for all $i \in I_{-1}$, and $Ad(b_k).\xi_{i,k}= e^{-2t_k}\zeta_{i,k}$ for all $i \in I_{-2}$. By the same calculation as before, we get that:

-  $||Ad(b_k).\Omega_{\hat x_k}(\xi_{\hx_k,i,k},\xi_{\hx_k,j,k})|| \leq e^{-2 t_k}\Omega_{\hat y_k}(\zeta_{\hy_k,i,k},\zeta_{\hy_k,j,k})$ for all $1\leq i <j \leq n$.

- $||Ad(b_k).\Omega_{\hat x_k}(\xi_{\hx_k,i,k},\xi_{\hx_k,j,k})|| \leq e^{-3 t_k}\Omega_{\hat y_k}(\zeta_{\hy_k,i,k},\zeta_{\hy_k,j,k})$ if $i$ or $j$ is in $I_{-2}$.  

Looking at the decomposition of $\lieg$ into  rootspaces, we see that  the converging sequences $(u_k)$ of $\lieg$ such that $||Ad(b_k).u_k||$  is contracted at a rate at least $e^{-2 t_k}$ must converge in $\lieh_{\infty}$, whereas  those sequences $(u_k)$ such that $||Ad(b_k).u_k||$  is contracted at a rate at least $e^{-3 t_k}$ must converge to $0_{\lieg}$. We infer that  $\Omega_{\hx_{\infty}}(\xi_{\hx_{\infty},i,\infty},\xi_{\hx_{\infty},j,\infty})=0$ as soon as $i$ or $j$ is in $I_{-2}$. If both $i$ and $j$ are in $I_{-1}$, we get that $\Omega_{\hx_{\infty}}(\xi_{\hx_{\infty},i,\infty},\xi_{\hx_{\infty},j,\infty})$ is in $l_{2, \infty}^{-1}.\liez^-$.  By the hypothesis of regularity on the connection $\omega$, this is possible only if $\Omega_{\hx_{\infty}}(\xi_{\hx_{\infty},i,\infty},\xi_{\hx_{\infty},j,\infty})=0$.  We finally get that $\Omega_{\hx_{\infty}}=0$ on  ${\lieh}_{\infty}$, and thus $\Omega_{\hx_{\infty}}=0$.     

\end{preuve}

Let us fix $x_{\infty}$ a point of the open set $U$ given by Proposition \ref{P}. Then, if $(x_k)$ is the  sequence constant to $x_{\infty}$, and $y_k=f_k(x_{\infty})$, then $(f_k,x_k,x_{\infty},y_{\infty})$ is a stable data. We note this peculiar stable data $(f_k,x_{\infty},x_{\infty},y_{\infty})$. 


\begin{proposition}
The point $x_{\infty}$ has an open neighbourhood $\Lambda$,  which is geometrically isomorphic to a  quotient $\Gamma \backslash \Omega_{o}$, with $\Gamma \subset P$ a dicrete subgroup (possibly trivial).

\end{proposition}

\begin{preuve}

We consider an holonomy data $(b_k,\hx_{\infty},\hx_{\infty},\hy{\infty})$ associated to $(f_k,x_{\infty},x_{\infty},y_{\infty})$. We still call $\phi_k=R_{b_k^{-1}} \circ {\hat f}_k$, and we use the same notations as in Lemma \ref{courbure}. In particular, the action of $(f_k)$ being equicontinuous at $x_{\infty}$, $b_k$ writes $l_{1,k}a^{t_k}l_{2,k}$ with $\lim_{k \to +\infty}t_k=+ \infty$.

 For every  $r \in \RR_+^*$, and  $k \in \NN \cup \{  \infty \}$, we denote by  $B_{r,k}$ the intersection of $\lieh_k$ with the ball of center  $0_{\lieg}$ and radius $r$  (for the norm  $||.||$). For every  $p \in B$, we set  $H_{p,{\lieh_k}}=\omega_p^{-1}(\lieh_k)$, and     $B_{p,r,k} = \omega_p^{-1}(B_{r,k})$. Since $\hat x_k$ and  $\hat y_k$ are both converging sequences, one can find $r \in \RR_+^*$ such that  $B_{r,k} \subset {\bf W}_{\hat x_k}  \cap {\bf W}_{\hat y_k} $ for every $k \in \NN \cup \{ \infty \}$.

\begin{lemme}
There is a sequence of integers $(k_m)_{m \in \NN}$ such that :

$\bullet$  for every  $m \in \NN$, $B_{\hx_{\infty},m,k_m} \subset W_{\hx_{\infty}}^B$, and $Exp_{\hx_{\infty}}$ maps $B_{\hx_{\infty},m,k_m}$ diffeomorphically on its  image (this image is denoted by  $U_m$).

$\bullet$ $V_m=[B_{m,k_m}] $ is a strictly increasing sequence of open subsets of $ {\bf X} \setminus \{  o^- \}$, the union of which is  ${\bf X} \setminus \{  o^- \}$

\end{lemme}

\begin{preuve}

From the relation $\phi_k \circ Exp_{\hx_{\infty}}=Exp_{\hy_k} \circ D_{\hx_{\infty}}\phi_k$, we infer that $D_{\hx_{\infty}}\phi_k$ maps $W_{\hx_{\infty}}^B$ on $W_{\hy_k}^B$. Lets us fix $m$. Since $Ad(b_k)(B_{m,k})$ tends to $0_{\lieg}$ as $k \to + \infty$,  $D_{\hx_{\infty}}\phi_k(B_{\hx_{\infty},m,k})$ tends to $0_{\hy_{\infty}}$ as $k \to + \infty$. Thus,  there exists an integer  $k_{1,m} \in \NN$ such that  $Exp_{\hy_k}$ maps $D_{\hx_{\infty}}\phi_k(B_{\hx_{\infty},m,k})$ diffeomorphically on its image, when  $k \geq k_{1,m}$. We infer that  $Exp_{\hx_{\infty}}$ maps $B_{\hx_{\infty},m,k}$ diffeomorphically on its image.  On the other hand, as  $k \to + \infty$, $[B_{m,k}]$ tends to $[B_{m, \infty}]$, which has compact closure in $\pi_{X} \circ Exp_G(Ad(\lieh_{\infty})={\bf X} \setminus \{  o^- \}$.  So, taking  $k \geq k_{2,m}$, we can suppose that $[B_{m,k}] \subset {\bf X} \setminus \{  o^- \}$.  Finally, since $[B_{m, \infty}]$ is a strictly increasing sequence, one can find a sequence $(k_m)$, with  $k_m \geq \max (k_{1,m}, k_{2,m})$ for every $m$, such that $[B_{m,k_m}]$ is also a strictly increasing sequence, for the inclusion. \end{preuve}

The "North-South" dynamics properties for the sequence  $(b_k)$ imply that   $\lim_{k \to + \infty} b_k.[B_{m,k_m}] = [o]$. Lemma \ref{collapse} then gives:
$\lim_{k \to + \infty} f_k([B_{\hx_{\infty},m,k_m}])=[y_{\infty}]$. We conclude that the action of $(f_k)$ is equicontinuous at each point $x \in U_m$. As a consequence of Proposition \ref{courbure},   the curvature vanishes on $U_m$.

We are now in the following situation.  The open set  $U_{\infty}=\bigcup_{m=1}^{\infty}U_m$ is flat and contains $x_{\infty}$. Let ${\tilde U}_{\infty}$ be its universal cover. Let  $\tilde x_{\infty} \in \tilde U_{\infty}$ be a point over  $x_{\infty}$, and  $\delta : \tilde U_{\infty} \to \X$ a developping map, mapping $\tilde x_{\infty}$ on $o$ (see section  \ref{plat}). For every $m \in \NN$, $U_m$ can be lifted to an open subset $\tilde U_m \subset \tilde U_{\infty}$ in the following way:  for every $\zeta \in B_{\hx_{\infty},m,k_m}$, the curve  $\zeta^*$ has a unique lift ${\tilde \zeta}^* \in C^1({\tilde U}_{\infty},[0,1],\tilde x_{\infty})$. The open set  $\tilde U_m$ is the union of such lifts. Since for every $\zeta \not = \xi$ in $B_{\hx_{\infty},m,k_m}$, $\zeta^*(1) \not = \xi^*(1)$ (indeed $Exp_{x_{\infty}}$ maps $B_{\hx_{\infty},m,k_m}$ diffeomorphically on its image), the projection of $\tilde U_m$ on $U_m$ is a diffeomorphism.

We use the results of section  section \ref{plat}. Let ${\hx}_{\infty}^{\prime}=R_{b_0}.{\hx}_{\infty}$ such that $\tilde \delta({\hx}_{\infty}^{\prime})=e$. Then    $\delta(\tilde \zeta^*)={\cD}_{x_{\infty}}^{{\hx}_{\infty}^{\prime}}(\zeta^*)=b_0.{\cD}_{x_{\infty}}^{{\hx}_{\infty}}(\zeta^*)$, so that $\delta(\tilde U_m)=b_0.V_m$. 
\begin{lemme}
The map  $\delta$ is injective on  $\tilde U_{\infty}$, and realizes a geometrical isomorphism between $\tilde U_{\infty}$ and $\X \setminus \{ b_0. o^- \}$.  
\end{lemme}

\begin{preuve}
We first prove that for all $m \in \NN$, $\delta$ is injective on $\tilde U_m$. Let $\tilde x$ and  $\tilde y$ be distinct points in  $\tilde U_m$. They project on two distinct points  $x$ and $y$ in $U_m$. We write $x=Exp_{x_{\infty}}(\zeta_{\hx_{\infty},m,k_m})$ and $y=Exp_{x_{\infty}}(\xi_{\hx_{\infty},m,k_m})$ with $\zeta \not = \xi$ in $B_{m,k_m}$. Thus, $\delta(x)= \zeta^*(1)$ and $\delta(y)=\xi^*(1)$. But $\pi_X \circ Exp_G$ is a diffeomorphism from $B_{m,k_m}$ on its image. As a consequence,  $\delta(x) \not = \delta(y)$.

To conclude, we just remark that  $\tilde U_{m}$ is an increasing sequence of open subsets.  \end{preuve}

It follows from the previous lemma that  $\tilde U_{\infty}$ is geometrically isomorphic to the open set  $\X \setminus \{  b_0.o^-\}$, which is itself geometrically isomorphic to $\Omega_o$. We get that  $U_{\infty}$ is geometrically isomorphic to a quotient $\Gamma \backslash \Omega_{o}$, where $\Gamma \subset P$ is a discrete subgroup acting freely properly discontinuously on $\Omega_o$.

\end{preuve}

We conclude the proof of Theorem \ref{theoreme1} thanks to Lemma  \ref{bebe-obata}, and the following result, applied to the inclusion $U_{\infty} \subset M$ :

\begin{theoreme}
\label{plongement}

Let $(M,B,\omega)$ be a Cartan geometry modelled on  $\X= \partial {\bf H}_{\KK}^d$. Let $\Gamma$ be a discrete subgroup of $P$, acting freely properly discontinuously on  $\Omega_o$. We suppose that there is a geometric embedding  $\sigma$ from $\Gamma \backslash \Omega_o$ into $M$.  Then we are in one of the following cases :

$(i)$ $M$ is geometrically isomorphic to the model space ${\bf X}$. In this case, $\Gamma = \{ e \}$, and there is a   point $\kappa \in M$ such that  $\sigma(\Omega_o)=M \setminus \{  \kappa \}$.

$(ii)$  The embedding $\sigma$ is a geometrical isomorphism between  $\Gamma \backslash \Omega_{o}$ and $M$.
\end{theoreme}

The proof of this result  will be the aim of the last section of this article. As an illustration, let us remark that Theorem  \ref{plongement}, when applied to the case of conformal riemannian structures, yields the:   

\begin{corollaire}
Let us suppose that $\sigma$ is a conformal embedding from a flat, complete, Riemannian manifold  $M$  of dimension $n \geq 3$,  into a Riemannian manifold $N$ of the same dimension. Then, either $\sigma$ is a conformal diffeomorphism between $M$ and $N$, or $M$ is the Euclidean space $E_n$, and  $N$  is conformally diffeomorphic to the standard sphere ${\bf S}^n$. In this case, $\sigma$ is the composition of the standard stereographic projection, with some M\"obius transformation.

\end{corollaire}

\section{A rigidity theorem for the geometric embeddings}
This last part is devoted to the proof of Theorem \ref{plongement}. Other cases of rigidity for geometrical embeddings of certain Cartan geometries will be studied more extensively  in \cite{frances}.

\subsection{Geometrical embeddings}

Let $(M,B,\omega)$ and $(N,B^{\prime},\omega^{\prime})$ be two manifolds endowed with Cartan geometries modelled on $\X=G/P$ (we are speaking here of the most  general framework : $G$ is any Lie group and $P$ is a closed subgroup of  $G$). In particular, $M$ and $N$ have the same dimension. 

A map  $\sigma$ will be called  {\it a geometrical embedding } from $M$ to $N$ if $\sigma$ lifts into a fiber bundle embedding $\hat \sigma : B \to B^{\prime}$, satisfying  $\hat \sigma^* \omega^{\prime}=\omega$.
 
\subsection{Cauchy boundary of a  Cartan geometry}

We decribe here a way for attaching a boundary to any Cartan geometry $(M, B, \omega)$, modelled on a space $\X=G/P$. This method was already used in general relativity, to associate a boundary to spacetimes (see for example  \cite{schmidt}).

We already saw in section \ref{cartan} that the choice of a basis of $\lieg$ yielded a parallelism  ${\cal R}$, and a Riemannian metric $\rho$ on $B$. The parallelism is orthonormal for the metric $\rho$. Let us denote by  $d_{\rho}$ the distance associated to the metric $\rho$. We can consider  $(\overline B, \overline d)$,  {\it the Cauchy completion } of the metric space $(B,d_{\rho})$. For every  $b \in P$, the differential of $R_b$, when expressed in the trivialisation of $TB$ given by the parallelism ${\cal R}$,  is just the linear transformation  $Ad(b^{-1})$. As a consequence,   $R_b$ is a uniformly continuous transformation of  $(B,d_{\rho})$ for all $b \in P$. Thus, it  can be prolongated into a continuous transformation    $\overline R_b$ of $(\overline B, \overline d)$. The space $\overline B$ writes naturally as the union $\overline B = B \cup \partial_c B$. The  {\it Cauchy boundary} $\partial_cB$ is stable for the action of the $\overline R_b$'s, $b \in P$. We then set $\partial_cM= \partial_c B/P$ (the quotient is taken for the action of  $P$ on  $\partial_c B$ by the transformations $\overline R_b$).

Of course, in general, the space $M \cup \partial_cM= \overline B/P$ behaves very badly from a topological point of view (it is generally non Hausdorff), so that the previous construction has a limited interest to get "nice boundaries".

Let us illustrate this construction in the case  $\X= \partial {\bf H}_{\KK}^d$.  

\begin{lemme}
\label{bord}
\ \\
\noindent $(i)$ $\partial_c \X= \emptyset$.

\noindent $(ii)$ If $\Gamma \subset G$ is a discrete subgroup acting freely properly discontinuously on  $\Omega_o$, the boundary $\partial_c ( \Gamma \backslash \Omega_o)$  reduces to a  point.

\end{lemme}

\begin{preuve}
Let us recall that for the model space  ${\bf X}$, the Cartan bundle is just the group  $G$, and the choice of a basis of  $\lieg$ determines a unique left invariant Riemannian metric $\rho_G$ (making the given basis an orthonormal basis). We denote by $B_o$ the preimage of ${\Omega_o}$ by the  fibration $\pi_X : G \to {\bf X}$. The topological boundary $\partial B_o$ of this open subset in $G$ is just the fiber over $o$.  Let  $\rho_o$ be the metric induced by   $\rho_G$ on $B_o$.

The Riemannian manifold $(G,\rho_G)$ is homogeneous, hence complete, what proves  point $(i)$. The group $\Gamma$ acts by isometries for $\rho_G$. The manifold  $\Gamma \backslash G$ is endowed with a metric $\overline \rho_G$, induced by  $\rho_G$ via the covering map. 
The manifold  $(\Gamma \backslash G,\overline \rho_G)$ is also complete. The  quotient $\overline B_o = \Gamma \backslash B_o$ can be identified with an open subset of $\Gamma \backslash G$, and its topological boundary $\partial \overline B_o \subset \Gamma \backslash G$ is a submanifold $\Gamma \backslash \partial B_o$. 
Since $\partial \overline B_o$ is a closed submanifold of $\Gamma \backslash G$, any point $p_{\infty} \in \partial \overline B_o $ is obtained as the limit of a Cauchy sequence  $(p_n)$ of $\overline B_o$ (for $d_{\overline \rho_o}$, the metric induced by $\rho_o$ via the covering map). 
Reciprocally, every sequence of  $\overline B_o$, which is a Cauchy sequence for  $d_{\overline \rho_o}$, is also a Cauchy  sequence for $d_{\rho_G}$. 
Thus, $(\Gamma \backslash G, d_{{\overline \rho}_G})$ is the Cauchy completion of   $(\overline B_o, d_{{\overline \rho}_o})$, and $\partial_c \overline B_o$ can be identified with  $\partial \overline B_o$. 
Since $P$ acts transitively on  $\partial \overline B_o$, we get that $\partial_c(\Gamma \backslash \Omega_o)$ is just one  point.

\end{preuve}

\subsection{The geometrical boundary of an embedding}

Let $(M, B, \omega)$ and  $(N, B^{\prime}, \omega^{\prime})$ be two Cartan geometries modelled on the same space  $\X=G/P$.

We assume that there is a geometrical embedding  $\sigma : M \to N$. Then, the topological boundary of $\sigma(M)$ in $N$ is called the  {\it geometrical boundary of $M$, associated with the embedding  $\sigma$}, and denoted $\partial_{\sigma}M$.  We will also write  $\partial_{\hat \sigma}B$ for the topological boundary of $\hat \sigma(B)$ in $B^{\prime}$. 

We call $\rho$ and  $\rho^{\prime}$ the Riemanniann metrics on $B$ and  $B^{\prime}$, associated  {\it to the choice of a same basis of  $\lieg$}, as they were defined in section \ref{metrique}. For these metrics, the embedding  $\hat \sigma$ is in fact an isometric embedding: $\hat \sigma^*\rho^{\prime}=\rho$. We call $d_{\rho}$ and  $d_{\rho^{\prime}}$ the distances associated to $\rho$
 and $\rho^{\prime}$.

 \begin{definition}
We call  $d_{\rho^{\prime}}^{\sigma}$ the distance on  $\hat \sigma(B)$ defined by:

\[  \dsigma(p,q) = \inf \{ L_{\rho^{\prime}}(\gamma) \ | \ \gamma \in C^1([0,1], \hat \sigma(B)), \ \gamma(0)=p, \ \gamma(1)=q \}  \] 
 
 \end{definition}
 
 Let us remark that  $\hat \sigma$ is an isometry from  $(B,d_{\rho})$ to $(\hat \sigma(B),\dsigma)$. In particular, by equivariance of the action of  $P$, the transformations $R_{b}$, for $b \in P$, are uniformously continuous for the distance $\dsigma$.

 We will need also the following:
 
 \begin{definition}[Regular points]
 We say that a point  $p \in \partial_{\hat \sigma}B$ is regular if there exist a sequence $(p_n)$ of $\hat \sigma(B)$ which tends to $p$, and such that $(p_n)$ is a Cauchy sequence for the distance $\dsigma$. 
 
 - A point $x \in \partial_{\sigma}M$ is said to be regular if there exist a regular point  $\partial_{\hat \sigma}B$ over $x$. 
 
 - The set of regular ponts of  $\partial_{\hat \sigma}B$ (resp. $\partial_{\sigma}M$) is denoted by $\partial_{\hat \sigma}^{reg}B$ (resp. $\partial_{\sigma}^{reg}M$). 
 
 \end{definition}
Let us remark that since  
 $\partial_{\hat \sigma}^{reg}B$ is invariant under the action of  $P$ on $B^{\prime}$, if there exist a point of   $\partial_{\hat \sigma}^{reg}B$ over $x$, then all the points of the fiber over  $x$ are regular.  
 
 \begin{lemme}
\label{regulier}
If $\partial_{\sigma}M$ is not empty, then   $\partial_{ \sigma}^{reg}M$ is dense in $\partial_{\sigma}M$.  
 
 \end{lemme} 

 \begin{preuve}
 Let us pick a  point $x_{\infty} \in \partial _{\sigma}M$, and lift it into a  point $\hx_{\infty} \in \partial_{\hat \sigma}B$. We fix a small ball $B_{\rho^{\prime}}(0_{\hx_{\infty}},\epsilon) \subset T_{\hx_{\infty}}B^{\prime}$ such that $Exp_{\hx_{\infty}}$ maps $B_{\rho^{\prime}}(0_{\hx_{\infty}},\epsilon) $ diffeomorphically on its  image. Then, there exists $u_0 \in B_{\rho^{\prime}}(0_{\hx_{\infty}},\epsilon) $ such that ${\hat u}_0^*(]0,1]) \cap \hat \sigma(B) \not = \emptyset$. Indeed, if it were not the case, the open set  $Exp_{x_{\infty}}(B_{\rho^{\prime}}(0_{\hx_{\infty}},\epsilon) ) $  would not intersect $ \sigma(M)$, contradicting  $x_{\infty} \in \partial_{\sigma}M$. 
 
Now, if there is a $t_0 \in ]0,1]$ such that  ${\hat u}_0^*(]0,t_0]) \subset \hat \sigma(B)$, it means that $\hx_{\infty}$ itself is regular. Indeed, if $(t_n)$  is a sequence of $]0,t_0]$ tending to  $0$, the sequence ${\hat u}_0^*(t_n)$ is a Cauchy sequence for $\dsigma$ and tends to $\hx_{\infty}$. 
 
 If a $t_0$ as above does not exist, then ${\hat u}_0^*(]0,1]) \cap \hat \sigma(B)$ has infinitely many connected components. There is a decreasing sequence of $]0,1]$, let us call it  $(t_n)_{n \in \NN}$,  converging to $0$, and such that those connected components are the intervals  ${\hat u}_0^*(]t_{2k+1},t_{2k}[)$, $k \geq 0$. By the same argument as above, the  points ${\hat u}_0^*(t_n)$ are in  $\partial_{\hat \sigma}^{reg}B$, for $n \geq 1$. We get a sequence of regular points tending to  $\hat x_{\infty}$. Projecting on $N$, we get a sequence of $\partial_{\sigma}^{reg}M$ converging to  $x_{\infty}$.
 
 \end{preuve}

 \begin{lemme}
\label{injection} 
There is  a map   $\hat {\j} : \partial_{\hat \sigma}^{reg}B \to {\cal P}(\partial_cB)$ (where  ${\cal P}(\partial_cB)$ stands for the set of parts of $\partial_cB$), such that if  $\hx$ and  $\hy$ are two distinct points of $\partial_{\sigma}^{reg}B$, then ${\hat {\j}}(x) \cap{\hat {\j}}(y) = \emptyset$. 

Moreover, the map  $\hat {\j}$ is equivariant for the action of  $P$ on $\partial_{\hat \sigma}^{reg}B $ and ${\cal P}(\partial_cB)$. It thus induces an injective  map $j : \partial_{\sigma}^{reg}M \to {\cal P}(\partial_cM)$.
 
 \end{lemme}

\begin{preuve}
We define the map  $\hat {\j}$  in the following way: for every  point $\hx_{\infty} \in \partial_{\hat \sigma}^{reg}B $, ${\hat {\j}}(\hx_{\infty})$ is the set of sequences $(\hx_k)$ of  $\hat \sigma(B)$, converging to  $\hx_{\infty}$,  and which are  Cauchy sequences for $\dsigma$ (to be precise, the points of    $\partial_cB$ are rather defined by the sequences $\hat \sigma^{-1}(\hx_k)$). With this definition, if $b \in P$, it is clear that $\hat {\j}(R_b.\hx_{\infty})$ is the part $R_b.\hat {\j}(\hx_{\infty})$, showing the equivariance of $\hat {\j}$. If  $\hx_{\infty}$ and  $\hy_{\infty}$  are distinct in $\partial_{\hat \sigma}^{reg}B$, and if $(\hx_k)$ and  $(\hy_k)$ are two sequences of  $\hat {\j}(\hx_{\infty})$ and  $\hat {\j} (\hy_{\infty})$ respectively, then there is $\epsilon > 0$ such that for  $k \in \NN$, $d_{\rho^{\prime}}(\hx_k,\hy_k)> \epsilon$. Hence, {\it a fortiori},  $\dsigma (\hx_k,\hy_k)> \epsilon$, what proves that  $(\hx_k) \not = (\hy_k)$ in $\partial_cB$. 

Let us now check that the induced map $j$ is also injective. Let us pick $x \not = y$ in $ \partial_{\sigma}^{reg}M$, and choose $\hat x$ and $\hy$ over $x$ and $y$ respectively. If $j(x) \cap j(y) \not = \emptyset$, it means that ${\hat {\j}}(\hx) \cap R_b({\hat {\j}}(\hy)) \not = \emptyset$, for some $b \in P$. But by equivariance, it means that ${\hat {\j}}(\hx) \cap {\hat {\j}}(R_b.\hy) \not =  \emptyset$, a contradiction since $\hx \not = R_b. \hy$. 

\end{preuve} 
 
\subsection{Proof of Theorem  \ref{plongement}}

We assume that the hypothesis of theorem  \ref{plongement} are satisfied, and we call $M$ the quotient $\Gamma \backslash \Omega_o$. 

The first case to deal with is  $\partial_{\sigma}M = \emptyset$. In this case, $\sigma$ is a diffeomorphism. It is thus a geometric isomorphism between   $M$ and $N$, and we are in the case  $(ii)$ of the theorem. 

Now, let us assume  $\partial_{\sigma}M \not = \emptyset$. By Lemma \ref{regulier}, the set of regular points of  $\partial_{\sigma}M$ is dense in  $\partial_{\sigma}M$, and in particular, is nonempty. But Lemma \ref{injection}, together with  point $(ii)$ of Lemma \ref{bord},  ensures that the set of regular points is a  singleton, and moreover, has to be dense in  $\partial_{\sigma}M$. This proves that  $\partial_{\sigma}M$ itself has just one point, that we call  $\kappa$. Thus, $\partial_{\hat \sigma}^{reg}B= \partial_{\hat \sigma}B$ can be identified with the fiber  (i.e a $P$-orbit)  of $B^{\prime}$ over  $\kappa$. The map $\hat {\j}$ of Lemma \ref{injection} is an equivariant map from this fiber, onto $\partial_c B=\Gamma \backslash P$. This forces  the action of $P$ on  $\Gamma \backslash P$  to be free, which implies $\Gamma = \{ e \}$.  

Finally, $N= \sigma(\Omega_o) \cup \{  \kappa \}$. So, $N$ is diffeomorphic to  $\X$, hence simply connected.  Moreover, $(M,B,\omega)$ is flat since it is flat on a dense open set. Thus, there is a developping map from  $N$ to $\X$, which turns out to be a covering map, by compacity of $N$. This developping map is then a geometrical isomorphism between $N$ and $\X$.

\ \\
Charles FRANCES\\
Laboratoire de Topologie et Dynamique\\
Universit\'e de Paris-Sud, B\^at. 430\\
91405 ORSAY\\
FRANCE\\
email: Charles.Frances@math.u-psud.fr\\
\end{document}